\documentclass[11pt,oneside,a4paper]{article}
\usepackage[T1]{fontenc}
\usepackage{lmodern}
\usepackage{xparse}
\usepackage{amsmath}
\usepackage{amssymb}
\usepackage[english]{babel}
\usepackage[utf8]{inputenc}
\usepackage[]{nicefrac}
\usepackage{amsmath, amssymb, amsthm}
\usepackage{graphicx}
\usepackage[colorinlistoftodos]{todonotes}
\usepackage{tikz}
\usepackage[affil-it]{authblk}
\usepackage[]{ccicons}
\usepackage[]{xcolor}
\usepackage{leftidx}
\usepackage{natbib}
\usepackage{pgfplots}
\usepackage[]{hyperref}
\usepackage{xcolor}
\hypersetup{%
    colorlinks = true,
    allcolors = gray,
}

\DeclareMathOperator*{\argmin}{argmin} 

\newcommand{\bxi}{\boldsymbol{\xi}}

\newcommand{\dVblank}{\, {\rm d}V}
\newcommand{\dxi}{\, {\rm d}\bxi}

\newcommand{\norm}[1]{\left\lVert#1\right\rVert}

\NewDocumentCommand \dV{ o }{
    \IfNoValueTF{#1}{\dVblank}
    {
        \dV_{#1}
    }
}

\NewDocumentCommand \mc{ m }{
    \mathcal{#1}
}

\NewDocumentCommand \eref{ m }{
    (\ref{eqn:#1})
}

\NewDocumentCommand \an{ m }{
    \langle {#1} \rangle
}

\NewDocumentCommand \mbf{ m }{
    \mathbf{#1}
}

\NewDocumentCommand \bs{ m }{
    \boldsymbol{#1}
}

\NewDocumentCommand \stateOne{ m }{
    \underline{\mbf{#1}}
}

\NewDocumentCommand \stateTwo{ m o }{
    \IfNoValueTF{#2}{\stateOne{#1}}
    {
    \stateOne{#1}\an{#2}
    }
}

\NewDocumentCommand \s{ m o o }{
    \IfNoValueTF{#3}{\stateTwo{#1}[#2]}
    {
    \underline{\mbf{#1}}(\mbf{#3})\an{#2}
    }
}

\NewDocumentCommand \stateOneI{ m m }{
    \underline{#1}_{#2}
}

\NewDocumentCommand \stateTwoI{ m m o }{
    \IfNoValueTF{#3}{\stateOneI{#1}{#2}}
    {
    \stateOneI{#1}{#2}\an{#3}
    }
}

\NewDocumentCommand \si{ m m o o }{
    \IfNoValueTF{#4}{\stateTwoI{#1}{#2}[#3]}
    {
        \stateOneI{#1}{#2}(\mbf{#4})\an{#3}
    }
}

\title{A machine-learning framework for peridynamic material models with physical constraints}
\date{\today}
\author{Xiao Xu\thanks{The Oden Institute for Computational Engineering and Sciences, The University of Texas at Austin, TX, xiaoxu42@utexas.edu}
\and Marta D'Elia \thanks{Computational Science and Analysis, Sandia National Laboratories, CA, mdelia@sandia.gov}
\and John T. Foster\thanks{The Oden Institute for Computational Engineering and Sciences, The University of Texas at Austin, TX, john.foster@utexas.edu}}

\begin{document}
\maketitle
\begin{abstract}
    As a nonlocal extension of continuum mechanics, peridynamics has been widely and effectively applied in different fields where discontinuities in the field variables arise from an initially continuous body. 
    An important component of the constitutive model in peridynamics is the \emph{influence function} which weights the contribution of all the interactions over a nonlocal region surrounding a point of interest. Recent work has shown that in solid mechanics the influence function has a strong relationship with the heterogeneity of a material's micro-structure. However, determining an accurate influence function analytically from a given micro-structure typically requires lengthy derivations and complex mathematical models. To avoid these complexities, the goal of this paper is to develop a data-driven regression algorithm to find the optimal bond-based peridynamic model to describe the macro-scale deformation of linear elastic medium with periodic heterogeneity. We generate macro-scale deformation training data by averaging over periodic micro-structure unit cells and add a physical energy constraint representing the homogenized elastic modulus of the micro-structure to the regression algorithm. We demonstrate this scheme for examples of one- and two-dimensional linear elastodynamics and show that the energy constraint improves the accuracy of the resulting peridynamic model.
\end{abstract}
\section{Introduction}
First proposed by~\cite{silling2000reformulation}, peridynamic (PD) mechanics models replace the spatial derivatives in the classical conservation of momentum equation with an integral functional to determine the net internal force density on material points during deformation. The integral equation provides convenience in modeling deformation problems with displacement discontinuities (e.g.\ cracks). The original PD formulation that is widely used in literature is the so-called \emph{bond-based peridynamic model}, which uses a pairwise force functional to describe the interaction between material particles \cite{zheng2020bond}.  A more general theory of peridynamics for solid mechanics, called \emph{state-based peridynamics}, was later proposed by \cite{silling2007psa}; however, the simpler bond-based peridynamic models still have a wide range of applications in brittle fracture \cite{huang2015improved, wang20183}. Additionally, bond-based models are useful for demonstrating advancements in computational techniques/implementations, homogenization theory, and/or numerical analysis where the more complicated state-based theory may be unnecessarily distracting from the central advancement of the research. 

The equation of motion in the bond-based PD theory is given by
\begin{align}
    \rho(\mbf x)\ddot{\mbf u}(\mbf x,t) = \int_{\mc H} \mbf{f} \left( \bs{\eta}, \bxi \right) \dxi + \mbf b(\mbf x,t),
    \label{eqn:peri}
\end{align}
where 
\begin{equation*}
\bs{\eta} := \mbf u(\mbf x + \bxi) - \mbf u(\mbf x),
\end{equation*}
$\rho$ is the mass density, $\mc H$ is a neighborhood of $\mbf x$ where integration of bond forces is carried out, $\mbf u$ is the displacement vector field, $\mbf b$ is a prescribed body force density field, and $\mbf f$ is the pairwise force function determining the force density that the particle at $\mbf x+\bxi$ exerts on the particle located by the position vector $\mbf x$. Regarding the constitutive model for the force function $\mbf f$, a scalar function is used to assign weights that scale the pair-wise force associated with each $\bxi$. For example, in~\cite{silling2005meshfree} the bond-based PD model for a linear elastic solid undergoing small deformation defined $\mathbf{f}$ as 
\begin{align}
    \mbf{f}\left( \bs{\eta}, \bxi \right) = c\left(\bxi\right)s\frac{\bxi}{\vert\bxi\vert},
    \label{eqn:const}
\end{align}
where $c$ is an \emph{influence function} which can be thought of as a spring constant in this setting that is, in general, unique for each $\bxi$. The bond stretch $s$ is defined as 
\begin{equation*}
s:=\frac{\vert\bxi + \bs\eta \vert}{\vert\bxi\vert}.
\end{equation*}
With a special choice of $c$, (\ref{eqn:peri}) and (\ref{eqn:const}) converge to the classical Cauchy linear momentum equation for a linear elastic solid with a Poisson ratio of $\nicefrac{1}{4}$ (in three dimensions) when the region of integration $\mc H$ is reduced to an infinitesimal volume~\cite{silling2008convergence}.

In a general setting, it is clear that the influence function plays an important role in determining the magnitude of the interaction between two material points and the overall effect of nonlocality of the material.  It is well-established  that the influence function is a key factor contributing to the behavior of peridynamic material models, especially in the case of wave dispersion~\cite{weckner2005effect} and fracture~\cite{seleson2011role}; however, there is still no general technique to systematically determine the influence function for different materials or to justify its choice for a given application. 

A few attempts have been made; for example, \cite{delgoshaie2015non} showed that the multi-scale~connectivity of natural pore networks can be used to extract nonlocal kernel functions for use in continuum nonlocal diffusion models. \cite{DElia2016ParamControl} used an optimal control technique to identify the nonlocal diffusivity parameter for both peridynamic and fractional nonlocal models.~\cite{sridhar2018general} presented a general multi-scale elastodynamic framework based on the Floquet-Bloch transform.~\cite{wildmandiscrete} used a prescribed  dispersion relation to derive the corresponding discrete micro-modulus function coefficients.  More recently,~\cite{Aksoylu2020} developed a selection criterion for nonlocal kernels derived from Taylor expansions that best approximate the classical linear dispersion relations for elastic solids.~\cite{You2020aaai} introduced a nonlocal-equation-constrained optimization algorithm to identify the optimal influence function for wave propagation through a one-dimensional heterogeneous bar.

In an explanation on the origin of nonlocality in solid materials,~\cite{silling2014origin} demonstrated that nonlocality can arise from the small-scale heterogeneity that is excluded through an implicit or explicit homogenization procedure.  Inspired by this idea, \cite{xu2020deriving}~focused on one-dimensional elastodynamics of a periodic heterogeneous bar and built a theoretical method to determine the PD influence function from the micro-structure. However, this theoretical method is effectively intractable for higher dimensional problems. Recently, \cite{you2020data} used machine learning to develop invertible nonlocal models from high-fidelity synthetic data while guaranteeing the well-posedness of the learned operator. Combining ideas from the previous two references, this work incorporates data-driven methods together with computational homogenization techniques to derive PD influence functions based on high-fidelity data associated to a specific material's micro-structure.

In this work, we aim to learn the discrete PD influence function from data, circumventing analytical complexity of conventional computational homogenization theories and avoiding costly fine-scale simulations that treat the heterogeneities in the micro-structure explicitly. The synthetic high-fidelity data set consists of microscopic displacement data of a given material obtained by solving the dynamic equations at micro-scale with a highly resolved micro-structural finite element model (FEM). This data set is then coarse-grained in each unit cell and it is used as training data set. The discrete PD influence function is the result of a data-driven regression algorithm; the corresponding nonlocal model serves as a surrogate model to describe the mechanics of the material at the macro-scale.  One of our main contributions is to demonstrate that the use of a physically-justified energy constraint imposed on the objective function reduces the amount of training data needed for a highly accurate homogenized model. This constraint enforces the Hill-Mandel macro-homogeneity condition under isotropic deformation. Our results show that this regression scheme is able to robustly learn the discrete PD influence function and to improve the accuracy of the predictions compared to standard choices of the PD influence function. In fact, standard nonlocal influence functions are nonnegative; while this guarantees well-posedness, it compromises the ability to accurately predict the displacement~\cite{weckner2011determination}. Instead, the class of kernels learned by our algorithm is allowed to be sign-changing.

The rest of the paper is organized as follows: \S\ref{chap:problem} outlines the general multi-scale framework for a linear elasticity problem with a periodic heterogeneous micro-structure and the learning algorithm for discrete micro-modulus functions for the macro-scale nonlocal model. In \S\ref{chap:1Delastic} and \S\ref{chap:2Delastic}, the algorithm is applied to one- and two-dimensional elastodynamics problems and the energy constraints for both cases are introduced. \S\ref{chap:numerical} presents testing results of the learned PD influence functions and illustrates the benefit of using the energy constraint. \S\ref{chap:conclusions} summarizes our contributions and provides future research guidelines.

\section{General linear elasticity problem with periodic heterogeneity}

\label{chap:problem}

Consider the linear elastodynamics of an open bounded domain $\Omega\subset\mathbb{R}^n$, with $n$ being the number of space dimension. Let the elastic medium be periodically heterogeneous with microstructural length scale being much smaller than the length scale of the domain $\Omega$
\begin{align}
    l \ll L,
    \label{eqn:lengthscale}
\end{align}
where $L$ and $l$ are used to denote the length scale of domain $\Omega$ and the one of its heterogeneous micro-structures respectively. Let $\{\mbf r_i\}^n_{i=1}$ be the set of lattice vectors that describe the micro-structures unit cell. Then $\{\mbf r_i\}^n_{i=1}$ form a basis for $\mathbb{R}^n$ and we can define the periodic lattice $\mc R$ as
\begin{align}
    \mc{R} = \{ \mbf{v} \in \mathbb{R}^n \, \vert \, \mbf v = \sum^n_{i=1} a_i \mbf r_i, a_i \in \mathbb{Z}\}, \notag
\end{align}
where $\mathbb Z$ is the set of all integers. Naturally, every single unit cell in $\Omega$ can be specified using a unique set of integers $\mbf a = (a_1,a_2, \ldots,a_n)$ as
\begin{align}
    \mc{T}(\mbf a) = \{ \mbf v \in \mathbb R^n \, \vert \, \mbf v = \sum^n_{i=1} b_i \mbf r_i, b_i \in [a_i,a_i+1] \}. \notag
\end{align}
We assume that, for $t\in[0,T]$, only small deformations occur within the periodic medium $\Omega$ under external loads and the displacement of material points in the domain is governed by the linear elasticity equations. Therefore, the displacement field $\mbf u$ satisfies the following generalized Hook's law
\begin{align}
    \bs \sigma (\mbf x,t) = \mathbb{C}(\mbf x) \colon \frac{1}{2}\left(\nabla \mbf u(\mbf x,t) + \left(\nabla \mbf u(\mbf x,t)\right)^\intercal\right) \qquad \mbf x \in \Omega \quad {\rm{and}} \quad 0 \leq t \leq T, \notag
\end{align}
where $\bs \sigma$ is the stress tensor field, the colon is the double contraction operator\footnote{i.e.\ $\mbf{A}:\mbf{B} = A_{ij}B_{ij}$}, $\nabla$ is the spatial gradient operator and $\mathbb{C}$ is the material's fourth-order elasticity tensor. Due to the periodic nature of the material, the stiffness tensor $\mathbb{C}$ must satisfy
\begin{align}
    \mathbb{C}(\mbf x + \mbf r) = \mathbb{C}(\mbf x) \qquad \mbf x, \mbf x + \mbf r \in \Omega \quad {\rm{and}} \quad \mbf r \in \mc{R}. \notag
\end{align}
Similarly, the same periodicity condition holds for material's density $\rho$
\begin{align}
    \rho(\mbf x + \mbf r) = \rho(\mbf x) \qquad \mbf x, \mbf x + \mbf r \in \Omega \quad {\rm{and}} \quad \mbf r \in \mc{R}. \notag
\end{align}
The equation of linear momentum conservation is
\begin{subequations}
    \begin{align}
        \rho(\mbf x) \ddot{\mbf u}(\mbf x,t) &= \nabla \cdot \bs \sigma(\mbf x,t) + \mbf b(\mbf x,t) \quad \forall \, \mbf x \in \Omega \quad {\rm{and}} \quad 0\leq t\leq T, \\
        \mc B\mbf u(\mbf x,t) &= q(\mbf x,t) \quad \qquad \qquad \qquad \forall \mbf x \in \partial\Omega \quad {\rm{and}} \quad 0\leq t\leq T,
    \end{align}
    \label{eqn:elastodynamic eqn}
\end{subequations}
where the superimposed double dot on $\mbf u$ indicates double differentiation in time, $\nabla \cdot$ is the spatial divergence operator, $\mbf b$ is the body force density field over the domain $\Omega$ and $\mc B$ denotes a boundary operator ($\mc B$ being identity map corresponds to a Dirichlet boundary condition and $\mc B = \partial/\partial n$ a Neumann boundary condition). Note that the displacement $\mbf u(\mbf x,t)$ in the above equations represents the accurate displacement of every material point inside the domain $\Omega$. Thus, $\mbf u(\mbf x,t)$ should be treated as the micro-scale displacement in the multi-scale modeling scheme. With the purpose of describing the system by means of a courser-scale model, we introduce the macro-scale displacement $\bar{\mbf u}(\mbf x, t)$ as the average of the micro-scale displacement inside each unit cell.

To be more specific, for each $\mbf x \in \Omega$, there exists a unique set of integers $\mbf a_{\mbf x} = (a_1,a_2, \ldots,a_n)$ such that $\mbf x \in \mc T(\mbf a_{\mbf x})$; we define the macro-displacement at $\mbf x$ as the average micro-displacement of the unit cell $\mc T(\mbf a_{\mbf x})$\begin{align}
    \bar{\mbf u}(\mbf x, t) = \frac {\int_{\mc T(\mbf a_{\mbf x})} \mbf u(\bxi, t) \dxi }{ \int_{\mc T(\mbf a_{\mbf x})}  \dxi }.
    \label{eqn:defmacrodisp}
\end{align}
Other macro-scale quantities are defined in the same way. In the rest of the paper we use the bar notation, i.e. $\bar{\mbf u}$, to represent the macro-scale average over the unit cell. We propose that the macro-displacement $\bar{\mbf u}(\mbf x,t)$ of the solutions to the elastodynamic equation \eref{elastodynamic eqn} with periodically oscillating material constants satisfies the PD bond-based model
\begin{subequations}
    \begin{align}
        \bar\rho(\mbf x) \ddot{\bar{\mbf u}}(\mbf x,t) =  \mc L_\omega[\bar{\mbf u}](\mbf x) + \bar{\mbf b}(\mbf x,t) \qquad \quad  \mbf x \in \Omega, 0\leq t\leq T, \\
        \mc B_{I}\bar{\mbf u}(\mbf x,t) = \bar q(\mbf x,t)   \qquad \qquad\qquad \quad \qquad  \mbf x \in \Omega_I, 0\leq t\leq T,
    \end{align}
    \label{eqn:PDeqn}
\end{subequations}
where $\mc B_I$ is the corresponding nonlocal interaction operator specifying a volume constraint in an appropriate nonlocal interaction domain $\Omega_I$. 
For the bond-based PD model the nonlocal operator $\mc L_\omega$ is defined as
\begin{align}
     \mc L_\omega[\bar{\mbf u}](\mbf x) = \int_{\mc H} \omega(\bxi)(\bar{\mbf u}(\mbf x + \bxi,t) - \bar{\mbf u}(\mbf x,t)) \dxi,
     \label{eqn:discretePD}
\end{align}
where $\omega$ is the kernel function that determines the nonlocality of the operator and $\mc H$ is a neighborhood of $\mbf x$. In the model of peridynamics, a characteristic length scale $\epsilon$ called the medium's horizon is often used to identify the neighborhood $\mc H$ as
\begin{align}
    \mc H = \{\mbf x + \bxi \vert \mbf x \in \Omega,\mbf x + \bxi \in \Omega, \vert\bxi\vert < \epsilon\}. \notag 
\end{align}
As a consequence, for the well-posedness of problem \eqref{eqn:PDeqn} we require $\Omega_I$ to be a layer, or collar, of thickness $\epsilon$ surrounding the domain, see, e.g., \cite{Du2012}.

Owing to the definition of our macro-displacement~\eref{defmacrodisp}, the macro-displacement $\bar{\mbf u}(\mbf x)$ is a constant function inside each unit cell, which leads us to naturally discretize the nonlocal equation \eref{PDeqn} as the summation over the unit cells inside $\mc H$
\begin{align}
    \bar\rho_i \ddot{\bar{\mbf u}}_i(t) = \sum_{j\in \mc H_n} \omega_{i,j}(\bar{\mbf u}_j(t) - \bar{\mbf u}_i(t)) + \bar{\mbf b}_i(t) \quad i \in \Omega_n, 0\leq t\leq T, \notag 
\end{align}
where $\Omega_n$ is the enumeration of the unit cells in $\Omega$ and $\mc H_n$ is the numbering of the cells in $\mc H$. The subscript $i,j$ are used to identify the unit cells and, in $\omega_{i,j}$, they indicate the interacting cells.
Since the length scale of the unit cell is much smaller than the domain, as highlighted in \eref{lengthscale}, we simply define $\omega_{i,j}$ as the product of the value of the influence function at the unit cell and the volume of the unit cell itself. With this choice, in the PD literature $\omega_{i,j}$ is often referred to as discrete micro-modulus function between unit cell $i$ and unit cell $j$. Explicit definitions of the discrete micro-modulus function are introduced later on in the paper.

\subsection{The regression algorithm}\label{sec:general-regression}
In order to learn the kernel function $\omega(\bxi)$, we assume that we are given high-fidelity solutions to the elastodynamic equations~\eref{elastodynamic eqn} at micro-scale, which can be either generated by numerically solving the equations using FEM with a mesh that is refined enough to describe the micro-scale displacement inside the unit cell, or directly obtained from experimental data with high accuracy. 

{The first step is to homogenize such solutions as described in \eref{defmacrodisp}. This results in} a collection of macro-scale data of displacement and acceleration at several time instants
\begin{align}
    \{ \bar{\mbf u}(\mbf x,t_i) , \ddot{\bar{\mbf u}}(\mbf x,t_i) \}_{i=1}^{N_t}. \notag 
\end{align}
We use these data to train our macro-scale PD equation \eref{PDeqn} and obtain an optimal surrogate model $\omega^*$ for the kernel function. This is achieved by solving an optimization problem of the following form
\begin{align}
    \omega^* = \argmin_\omega \sum_{i=1}^{N_t} \norm{\mc L_\omega[\bar{\mbf u}(\mbf x,t_i)] + \bar{\mbf b}(\mbf x,t_i)-\bar{\rho}(\mbf x)\ddot{\bar{\mbf u}}(\mbf x,t_i)}_\Omega,
    \label{eqn:minimizationwocon}
\end{align}
where $\norm{\cdot}_\Omega$ denotes an appropriate norm over $\Omega$. Note that the above training procedure only makes use of the macro-scale dynamic data which can only describe the deformation of the medium under certain loads and boundary conditions; more importantly, the minimization problem does not contain any physical information about the micro-structure of the medium. Thus, the accuracy of the predictions corresponding to the optimal kernel highly depends on the type and number of training data. {This fact can compromise the generalization properties of this algorithm and, hence, the quality of the predictions. To overcome this limitation, we add physics-based information} to the cost function that take into account features of the micro-structure and, possibly, acts as regularizers for the optimization problem. With this addition, we reformulate \eqref{eqn:minimizationwocon} as 
\begin{align}
    \omega^* = \argmin_\omega \sum_{i=1}^{N_t} \norm{\mc L_\omega[\bar{\mbf u}(\mbf x,t_i)] + \bar{\mbf b}(\mbf x,t_i)-\bar{\rho}(\mbf x)\ddot{\bar{\mbf u}}(\mbf x,t_i)}_\Omega + \sum_{i=1}^{N_c} \alpha_i C_i(\omega),
    \label{eqn:minimizationwcon}
\end{align}
where $\alpha_i$ is the regularization parameter, $C_i$ is the corresponding physics-based constraint and $N_c$ is the number of constraints. In the next sections, the one- and two-dimensional elastodynamics models are further specified, and we present specific constraints, such as \eref{1Dconstraint} and \eref{2Dconstraint}, that describe different features of the micro-structure of the given medium.

\section{One-dimensional Elastodynamics}
\label{chap:1Delastic}

\begin{figure}
\centering
\begin{tikzpicture}[scale=0.65]
    \tikzstyle{pt}=[circle, fill=black, inner sep=0pt, minimum size=5pt]
    \draw (-8,0) -- (8,0) -- (8,0.8) -- (-8,0.8);
    \draw (-8,-1) -- (-8,1.8);
    \draw (-8,1.5) -- (-8.5,1);
    \draw (-8,1.2) -- (-8.5,0.7);
    \draw (-8,0.9) -- (-8.5,0.4);
    \draw (-8,0.6) -- (-8.5,0.1);
    \draw (-8,0.3) -- (-8.5,-0.2);
    \draw (-8,0) -- (-8.5,-0.5);
    \draw (-8,-0.3) -- (-8.5,-0.8);
    
    \draw(-3,0) -- (-3,0.8);
    \draw(-3.8,0) -- (-3.8,0.8);
    
    \draw[->] (-3.4,0.4)  to [out=220,in=160, looseness=1] (-4.2,-2.8);

	\draw[-latex] (-8,-0.4) -- (-6,-0.4); 		  
    \draw (-7.5,-0.7) node() {$x$};
    \draw[-latex] (8,0.4) -- (10,0.4) node[midway, above , fill=white, inner sep=1pt] {$load$} ;
    
    
   \filldraw[fill=black!50, draw=black] (-4,-2.5) rectangle (-3,-3.3);
    \filldraw[fill=black!50, draw=black] (-1,-2.5) rectangle (0,-3.3);
    \draw (-3,-3.3) -- (-1,-3.3);
    \draw (-3,-2.5) -- (-1,-2.5);
    
    \draw[latex-latex] (-4,-2.2) -- (-3,-2.2) node[midway , above,  fill=white, inner sep=1pt] {$l_1/4$} ;
    \draw (-4,-2.5) -- (-4,-1.9);
    \draw (-3,-2.5) -- (-3,-1.9);
    
    \draw[latex-latex] (0,-2.2) -- (-1,-2.2) node[midway , above,  fill=white, inner sep=1pt] {$l_1/4$} ;
    \draw (0,-2.5) -- (0,-1.9);
    \draw (-1,-2.5) -- (-1,-1.9);
    
    \draw[latex-latex] (-3,-2.2) -- (-1,-2.2) node[midway , above,  fill=white, inner sep=1pt] {$l_1/2$} ;
    
    \filldraw[fill=black!50, draw=black] (2,-2) rectangle (2.8,-2.8);
    \draw (4,-2.4) node() {$ E_s,\rho$};

    \filldraw[fill=white, draw=black] (2,-3) rectangle (2.8,-3.8);
    \draw (4,-3.4) node() {$ E_c,\rho$};

    

\end{tikzpicture}
\caption{One-dimensional composite with periodic heterogeneity}
\label{fig:1DIllustration}
\end{figure}
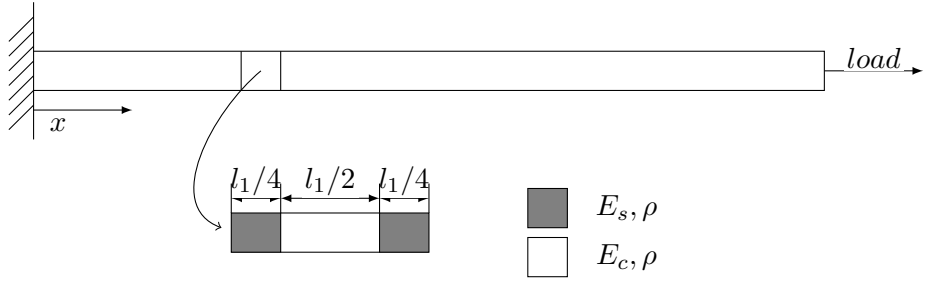

We consider a one-dimensional elastic composite problem similar to the one discussed in \citep{xu2020deriving}: a composite rod made of a periodic array of two linearly elastic, homogeneous, and isotropic constituents with perfect interfaces. The rod is fixed on one end and external loads are applied on the other end; {we also assume that no body forces are applied along the bar, i.e. $\mbf{b=0}$.} We define the micro-structure of the medium as the symmetric heterogeneous unit cell shown in Figure~\ref{fig:1DIllustration}, where $l_1$ is the length of the unit cell, the dark block represents the stiffer constituent with elastic modulus $E_s$ and density $\rho$ and the white block represents the more compliant constituent with elastic modulus $E_c$ and density $\rho$. If we use natural ordering from left to right as the numbering of the unit cells of the composite rod, we can simplify the discrete PD equation $\eref{discretePD}$ as 
\begin{align}\label{eq:discrete-1d}
    \rho\ddot{\bar u}_i(t) = \sum_{j=-\epsilon/l_1}^{\epsilon/l_1} \omega_j (\bar{u}_{i+j}(t) - \bar{u}_i(t)).  
\end{align}
Because of the symmetry of the medium's micro-structure, we enforce the symmetry constraint on the discrete micro-modulus function as follows:
\begin{align}\label{eq:symmetric-w}
    \omega_{-i} = \omega_i \qquad 0<i\leq\epsilon/l_1. 
\end{align}
By using the discrete expression in \eqref{eq:discrete-1d} in combination with \eqref{eq:symmetric-w}, the minimization problem \eref{minimizationwocon} reduces to the following linear regression problem for the vector $\mbf v_\omega := (\omega_{1},\omega_{2},\ldots,\omega_{\epsilon/l})$,
\begin{align}
    \mbf v_\omega^* = \argmin_{\mbf v_\omega}\sum_{i=1}^{N_t}\left(\sum_{j=1}^{\epsilon/l_1}\omega_j \left(\bar{u}_{i+j}(t_i)+\bar{u}_{i-j}(t_i) - \bar{u}_i(t_i)\right)-\ddot{\bar u}_i(t_i) \right)^2,
    \label{eqn:1Dminimization}
\end{align}
where we used the $\ell^2$ norm with respect to the values of the operator and the acceleration in each unit cell.

We next describe how to obtain a physics-based constraint that embeds the micro-structure information in the regression problem. In \citep{silling2005meshfree} the authors consider a large homogeneous body under isotropic extension and derive a relationship between the PD kernel function and the classical bulk modulus. Inspired by this work, we also assume that the composite rod is under constant strain, i.e. the macro-scale displacement can be expressed as $\bar u(x) = sx$, where $s$ is the corresponding macro-scale constant strain. As a consequence, the micro-scale displacement inside each constituent is also linear, which implies that the macro-scale constitutive equation is given by
\begin{align}
    \bar{\sigma} = E_{\rm{hom}} s.
    \label{eqn:homEs}
\end{align}
Here, $E_{\rm{hom}}$ is the homogenized elastic modulus of the unit cell
\begin{align}
    E_{\rm{hom}} = \frac{2}{1/E_s+1/E_c}. \notag 
\end{align}
In the framework of bond-based peridynamics, we can express the stress $\bar\sigma$ at point $x=0$ as the sum of all the bond forces in bonds that cross $x=0$ \citep{silling2003deformation}, i.e.
\begin{align}
    \bar{\sigma} &= \int^\infty_0\int^\infty_0 \omega(r+s)(\bar u(r) - \bar u(-s)) \rm ds \rm dr , \notag \\
    &= s \int^\infty_0 \xi^2 \omega(\xi) \rm d \xi.  
    \label{eqn:peridynamicstress}
\end{align}
Combining \eref{homEs} with \eref{peridynamicstress} and substituting discrete micro-modulus gives
\begin{align}
    E_{\rm{hom}} = \sum_{i=1}^{\epsilon/l_1} (i\;l_1)^2\omega_i. \notag 
\end{align}
Therefore, we define the constraint as
\begin{align}
    C(\mbf v_\omega) = \left( E_{\rm{hom}} - \sum_{i=1}^{\epsilon/l_1} (il_1)^2\omega_i \right)^2,
    \label{eqn:1Dconstraint}
\end{align}
and we refer to it as an {\it energy constraint}. Note that this corresponds to enforcing the effective elastic modulus of medium's micro-structures.

\section{Two-dimensional Elastodynamics}
\label{chap:2Delastic}
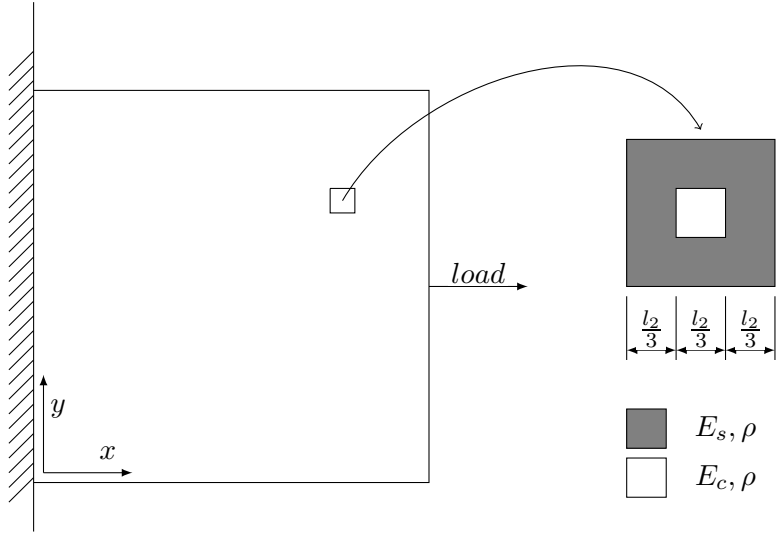
\begin{figure}
\centering
\begin{tikzpicture}[scale=0.65]
    \tikzstyle{pt}=[circle, fill=black, inner sep=0pt, minimum size=5pt]
    \draw (-8,-4) -- (0,-4) -- (0,4) -- (-8,4);
    \draw (-2,2) -- (-1.5,2) -- (-1.5,1.5) -- (-2,1.5) -- (-2,2);
    
    \draw (-8,-5) -- (-8,5.8);
    \draw (-8,4.8) -- (-8.5,4.3);
    \draw (-8,4.5) -- (-8.5,4.0);
    \draw (-8,4.2) -- (-8.5,3.7);
    \draw (-8,3.9) -- (-8.5,3.4);
    \draw (-8,3.6) -- (-8.5,3.1);
    \draw (-8,3.3) -- (-8.5,2.8);
    \draw (-8,3.0) -- (-8.5,2.5);
    \draw (-8,2.7) -- (-8.5,2.2);
    \draw (-8,2.4) -- (-8.5,1.9);
    \draw (-8,2.1) -- (-8.5,1.6);
    \draw (-8,1.8) -- (-8.5,1.3);
    \draw (-8,1.5) -- (-8.5,1);
    \draw (-8,1.2) -- (-8.5,0.7);
    \draw (-8,0.9) -- (-8.5,0.4);
    \draw (-8,0.6) -- (-8.5,0.1);
    \draw (-8,0.3) -- (-8.5,-0.2);
    \draw (-8,0) -- (-8.5,-0.5);
    \draw (-8,-0.3) -- (-8.5,-0.8);
    \draw (-8,-0.6) -- (-8.5,-1.1);
    \draw (-8,-0.9) -- (-8.5,-1.4);
    \draw (-8,-1.2) -- (-8.5,-1.7);
    \draw (-8,-1.5) -- (-8.5,-2.0);
    \draw (-8,-1.8) -- (-8.5,-2.3);
    \draw (-8,-2.1) -- (-8.5,-2.6);
    \draw (-8,-2.4) -- (-8.5,-2.9);
    \draw (-8,-2.7) -- (-8.5,-3.2);
    \draw (-8,-3.0) -- (-8.5,-3.5);
    \draw (-8,-3.3) -- (-8.5,-3.8);
    \draw (-8,-3.6) -- (-8.5,-4.1);
    \draw (-8,-3.9) -- (-8.5,-4.4);

    \draw[-latex] (0,0) -- (2,0) node[midway, above , fill=white, inner sep=1pt] {$load$} ;
    
    
    
    \filldraw[fill=black!50, draw=black] (4,3) rectangle (7,0);
    \filldraw[fill=white!50, draw=black] (5,1) rectangle (6,2);
    
    
    
    
    \draw[-latex] (-7.8,-3.8) -- (-6,-3.8);		  
    \draw (-6.5,-3.4) node() {$x$};
    
    \draw[-latex] (-7.8,-3.8) -- (-7.8,-1.8); 		  
    \draw (-7.5,-2.5) node() {$y$};
    
    \draw[latex-latex] (4,-1.3) -- (5,-1.3) node[midway , above,  fill=white, inner sep=1pt] {$\frac{l_2}{3}$};
    \draw (4,-1.5) -- (4,-0.2);
    \draw (5,-1.5) -- (5,-0.2);
    
    \draw[latex-latex] (6,-1.3) -- (7,-1.3) node[midway , above,  fill=white, inner sep=1pt] {$\frac{l_2}{3}$};
    \draw (6,-1.5) -- (6,-0.2);
    \draw (7,-1.5) -- (7,-0.2);
    
    \draw[latex-latex] (5,-1.3) -- (6,-1.3) node[midway , above,  fill=white, inner sep=1pt] {$\frac{l_2}{3}$};
    
    \filldraw[fill=black!50, draw=black] (4,-1-1.5) rectangle (4.8,-1.8-1.5);
    \draw (6,-1.4-1.5) node() {$ E_s,\rho$};
    
    \filldraw[fill=white, draw=black] (4,-2-1.5) rectangle (4.8,-2.8-1.5);
    \draw (6,-2.4-1.5) node() {$ E_c,\rho$};
    
    \draw[->] (-1.75,1.75)  to [out=60,in=120, looseness=1] (5.5,3.2);

    

\end{tikzpicture}
\caption{Two-dimensional composite with periodic heterogeneity}
\label{fig:2DIllustration}
\end{figure}

We extend the previous one-dimensional learning framework to two dimensions. We consider a square composite thin plate with one side fixed and external loads applied on the other side as shown in Figure~\ref{fig:2DIllustration}. {Also in this case, we assume that no external forces are applied, i.e. ${\bf b}={\bf 0}$.} The composite thin plate has periodic micro-structure whose lattice vectors are the unit vectors along the $x$- and $y$-axis. The micro-structure unit cell of the thin plate is a square composed of two different isotropic constituents; the compliant inclusions (represented by white block in Figure~\ref{fig:2DIllustration}) is embedded in the stiffer continuous matrix phase (represented by dark grey block in Figure~\ref{fig:2DIllustration}); the inclusions also have a square shape and are placed in the middle of the unit cell. Limited by the fact that the bond-based PD model can only be used to describe plane stress deformation of medium whose Poisson ratio is $\nicefrac{1}{3}$ \citep{trageser2020bond}, we assume that the Poisson ratio of every constituent in the medium is~$\nicefrac{1}{3}$ for consistency.
Following similar procedures as in {\S\ref{chap:1Delastic}}, the discretization of the two-dimensional PD equation for macro-scale displacements, the minimization problem \eref{minimizationwocon} can also be reduced to a linear regression problem for the discrete micro-modulus as \eref{1Dminimization}.  We choose the square neighborhood $\mc H$ of point $(x,y)$ 
\begin{align}
    \mc H = \{(x+\xi,y+\eta)\vert -\epsilon\leq\xi\leq\epsilon, -\epsilon\leq\eta\leq\epsilon \}. \notag 
\end{align}
{This choice is motivated by the specific micro-structure considered in this section and, more importantly, it does not compromise the well-posedness of the problem nor the convergence to the local limit as the nonlocal neighborhood shrinks, see, e.g. \cite{DElia2020cookbook}.}
We use natural ordering in $x$ direction and $y$ direction as the two-dimensional numbering of unit cells to discretize the PD equation for macro-scale displacements \eref{discretePD} as
\begin{align}
    \rho\ddot{\bar{\mbf u}}_{p,q}(t) = \sum_{i,j=-\epsilon/l_2}^{\epsilon/l_2} \omega_{i,j} (\bar{\mbf u}_{p+i,q+j}(t) - \bar{\mbf u}_{p,q}(t)). \notag 
\end{align}
Since the structure of the unit cell is symmetric about $i=0$, $j=0$ and $i=j$ and each constituent is isotropic, we can enforce the following symmetry constraints on our PD kernel function
\begin{align}
    \omega_{i,j} = \omega_{-i,j} = \omega_{i,-j} = \omega_{j,i} \qquad -\epsilon/l_2 \leq i,j\leq \epsilon/l_2. \notag 
\end{align}
Then, following similar procedures as in {\S\ref{chap:1Delastic}}, the minimization problem \eref{minimizationwocon} can also be reduced to a linear regression problem for the discrete micro-modulus function $\{\omega_{i,j} \}_{i,j=-\epsilon/l_2}^{i,j=\epsilon/l_2}$.

For the derivation of a two-dimensional energy constraint, we again consider the plane stress problem of the medium under isotropic extension, i.e. $\bar{\mbf u}(\mbf x) = s\mbf x$ where $s$ is a constant and compute the average strain energy density of the unit cell $W_{uc}$ from micro-scale displacement solutions (which, in our case, are obtained via FEM simulations). Once again, in the PD framework, by using the definition of the micro-potential \citep{silling2005meshfree}, the strain energy density can be expressed as 
\begin{align}
    W_{uc} = \frac{1}{2} \int_{\mc H} \frac{\omega(\bxi) s^2 \vert\bxi\vert^2}{2} \dxi. \notag 
\end{align}
Therefore, we define the energy constraint in two-dimensions as
\begin{align}
    C(\omega) = \left(W_{uc} - \frac{1}{2} \int_{\mc H} \frac{\omega(\bxi) s^2 \vert\bxi\vert^2}{2} \dxi \right)^2, 
    \label{eqn:2Dconstraint}
\end{align}
which weakly prescribes the value of the average strain energy density of the medium under isotropic extension (i.e.\ the effective bulk modulus of the micro-structure).

\section{Numerical Examples}
\label{chap:numerical}

In this section we use several one- and two-dimensional numerical examples to illustrate the performance of our PD kernel learning procedure and demonstrate the effectiveness of the energy constraint.

\subsection{One-dimensional test cases}
We consider the one-dimensional problem described in {\S\ref{chap:1Delastic}} with geometry parameters $L = 1 \textrm{m}$ and $l = 0.02 \textrm{m}$, and material properties $E_s = 200 \textrm{GPa}$, $E_c = 5 \textrm{GPa}$ and $\rho = 8000 \textrm{kg/m}^3$. We perform FEM simulations to calculate solutions to the small-scale elastodynamic equations \eref{elastodynamic eqn} with the following time-dependent boundary condition 
$$u_{bc}(t) = u_0 a_0 t^6(t-T_s)^6[1-H(t-T_s)],$$
where $u_0 = 1\times10^{-2} \mathrm{m}$, $a_0$ is a scaling factor, $H$ is the Heaviside function and $T_s = 0.157 \textrm{ms}$. We coarse-grain the micro-scale data into macro-scale data of displacement and acceleration at a set of discrete points in time over the interval $[0,T]$, for $T = 10^{-3}\rm s$. We divide this data set into a training set, for $t\in[0,T_t]$, and a testing data set, for $t\in(T_t,T)$ and we set $T_t = 0.17\times10^{-3}\rm s$. We choose the horizon of the PD model to be $\epsilon = 8l$; in order to avoid nonlocal boundary effects we discard the macro-scale data of points whose neighborhood is not fully contained in the domain.

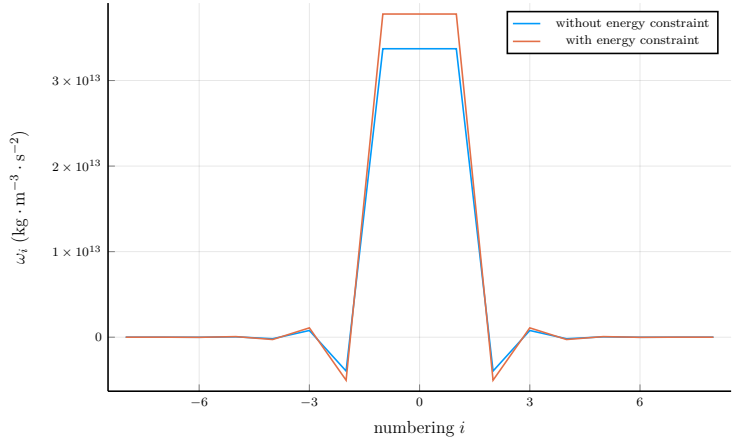
\begin{figure}[htbp]
\centering
\scalebox{0.6}{\begin{tikzpicture}[]
\begin{axis}[
  height = {101.6mm},
  ylabel = {$\omega_i\textrm{ }(\textrm{kg}\cdot \textrm{m}^{-3}\cdot\textrm{s}^{-2})$\\\color{white}.\\\color{white}.\\\color{white}.},
  xmin = {-8.48},
  xmax = {8.48},
  ymax = {3.905727762305043e13},
  xlabel = {$\textrm{numbering } i $},
  unbounded coords=jump,scaled x ticks = false,xlabel style = {font = {\fontsize{11 pt}{14.3 pt}\selectfont}, color = {rgb,1:red,0.00000000;green,0.00000000;blue,0.00000000}, draw opacity = 1.0, rotate = 0.0},xmajorgrids = true,xtick = {-6.0,-3.0,0.0,3.0,6.0},xticklabels = {$-6$,$-3$,$0$,$3$,$6$},xtick align = inside,xticklabel style = {font = {\fontsize{8 pt}{10.4 pt}\selectfont}, color = {rgb,1:red,0.00000000;green,0.00000000;blue,0.00000000}, draw opacity = 1.0, rotate = 0.0},x grid style = {color = {rgb,1:red,0.00000000;green,0.00000000;blue,0.00000000},
draw opacity = 0.1,
line width = 0.5,
solid},axis x line* = left,x axis line style = {color = {rgb,1:red,0.00000000;green,0.00000000;blue,0.00000000},
draw opacity = 1.0,
line width = 1,
solid},scaled y ticks = false,ylabel style = {align = center,font = {\fontsize{11 pt}{14.3 pt}\selectfont}, color = {rgb,1:red,0.00000000;green,0.00000000;blue,0.00000000}, draw opacity = 1.0, rotate = 0.0},ymajorgrids = true,ytick = {0.0,1.0e13,2.0e13,3.0e13},yticklabels = {$0$,$1\times10^{13}$,$2\times10^{13}$,$3\times10^{13}$},ytick align = inside,yticklabel style = {font = {\fontsize{8 pt}{10.4 pt}\selectfont}, color = {rgb,1:red,0.00000000;green,0.00000000;blue,0.00000000}, draw opacity = 1.0, rotate = 0.0},y grid style = {color = {rgb,1:red,0.00000000;green,0.00000000;blue,0.00000000},
draw opacity = 0.1,
line width = 0.5,
solid},axis y line* = left,y axis line style = {color = {rgb,1:red,0.00000000;green,0.00000000;blue,0.00000000},
draw opacity = 1.0,
line width = 1,
solid},    xshift = 0.0mm,
    yshift = 0.0mm,
    axis background/.style={fill={rgb,1:red,1.00000000;green,1.00000000;blue,1.00000000}}
,legend style = {color = {rgb,1:red,0.00000000;green,0.00000000;blue,0.00000000},
draw opacity = 1.0,
line width = 1,
solid,fill = {rgb,1:red,1.00000000;green,1.00000000;blue,1.00000000},fill opacity = 1.0,text opacity = 1.0,font = {\fontsize{8 pt}{10.4 pt}\selectfont}},colorbar style={title=},
  ymin = {-6.31307041653858e12},
  width = {152.4mm}
]

\addplot+[
  color = {rgb,1:red,0.00000000;green,0.60560316;blue,0.97868012},
draw opacity = 1.0,
line width = 1,
solid,mark = none,
mark size = 2.0,
mark options = {
            color = {rgb,1:red,0.00000000;green,0.00000000;blue,0.00000000}, draw opacity = 1.0,
            fill = {rgb,1:red,0.00000000;green,0.60560316;blue,0.97868012}, fill opacity = 1.0,
            line width = 1,
            rotate = 0,
            solid
        }
] coordinates {
  (-8.0, -2.138777855405617e8)
  (-7.0, 1.839021382397023e9)
  (-6.0, -9.870275055577085e9)
  (-5.0, 4.388432230080687e10)
  (-4.0, -1.8274343703620914e11)
  (-3.0, 7.829954421189196e11)
  (-2.0, -3.9555191939349336e12)
  (-1.0, 3.3716193774615156e13)
  (0.0, 3.3716193774615156e13)
  (1.0, 3.3716193774615156e13)
  (2.0, -3.9555191939349336e12)
  (3.0, 7.829954421189196e11)
  (4.0, -1.8274343703620914e11)
  (5.0, 4.388432230080687e10)
  (6.0, -9.870275055577085e9)
  (7.0, 1.839021382397023e9)
  (8.0, -2.138777855405617e8)
};

\addplot+[
  color = {rgb,1:red,0.88887350;green,0.43564919;blue,0.27812294},
draw opacity = 1.0,
line width = 1,
solid,mark = none,
mark size = 2.0,
mark options = {
            color = {rgb,1:red,0.00000000;green,0.00000000;blue,0.00000000}, draw opacity = 1.0,
            fill = {rgb,1:red,0.88887350;green,0.43564919;blue,0.27812294}, fill opacity = 1.0,
            line width = 1,
            rotate = 0,
            solid
        }
] coordinates {
  (-8.0, -3.713781520078972e8)
  (-7.0, 3.09907778033768e9)
  (-6.0, -1.6063129627098362e10)
  (-5.0, 6.855018816768585e10)
  (-4.0, -2.7157558915720724e11)
  (-3.0, 1.090520973175936e12)
  (-2.0, -5.029003962587948e12)
  (-1.0, 3.77732111690998e13)
  (0.0, 3.77732111690998e13)
  (1.0, 3.77732111690998e13)
  (2.0, -5.029003962587948e12)
  (3.0, 1.090520973175936e12)
  (4.0, -2.7157558915720724e11)
  (5.0, 6.855018816768585e10)
  (6.0, -1.6063129627098362e10)
  (7.0, 3.09907778033768e9)
  (8.0, -3.713781520078972e8)
};

\legend{{}{without energy constraint}, {}{with energy constraint}}
\end{axis}

\end{tikzpicture}}
\caption{Discrete micro-modulus function $\omega_i$ for the one-dimensional test case for both the unconstrained and constrained formulations.}
\label{fig:micro-modulus_1D}
\end{figure}

To test the efficacy of the energy constraint, we solve the linear regression problems both without and with the energy constraint, i.e.\ problems \eref{minimizationwocon} and \eref{minimizationwcon}, respectively. The optimal discrete micro-modulus functions $\{\omega^*_i\}_{i=-8}^{i=8}$ are reported in Figure~\ref{fig:micro-modulus_1D}. It is important to note that these two discrete micro-moduli generate positive definite matrices when discretizing the PD equations \eref{PDeqn}; this ensures the well-posedness of solutions to \eref{PDeqn}. 

\begin{figure}[htbp]
\centering
\scalebox{0.65}{\input{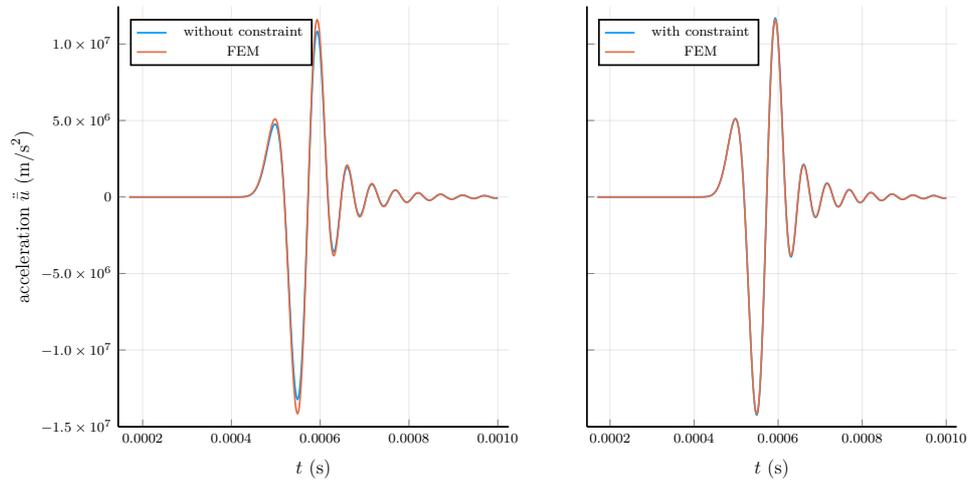}}
\caption{Testing results on the acceleration of the middle unit cell}
\label{fig:test_mid_acc_1D}
\end{figure}

\begin{figure}[htbp]
\centering
\scalebox{0.65}{\input{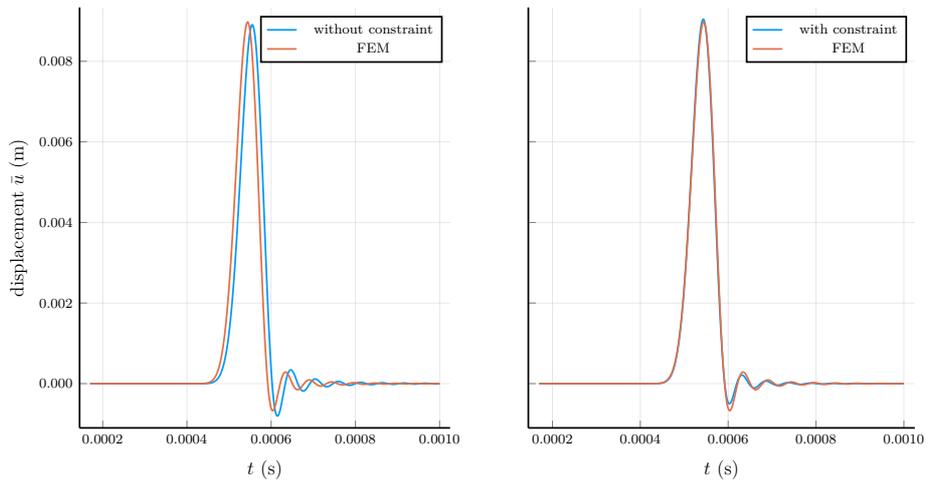}}
\caption{Testing results on the displacement of the middle unit cell}
\label{fig:test_mid_disp_1D}
\end{figure}

To test the performance of the learning algorithm, we use the two discrete micro-moduli to predict the macro-scale acceleration for $t \in (T_t,T]$. In Figure~\ref{fig:test_mid_acc_1D} we report the predicted acceleration of the unit cell at the middle of the rod as a function of time; as a reference, we also report the macro-scale result directly calculated from FEM simulations. 
The predicted acceleration using the micro-modulus with the energy constraint shows good agreement with the FEM solutions while the prediction using micro-modulus without the energy constraint has small but noticeable differences. We further test these two micro-moduli by predicting the deformation of the bar after time $T_t$, i.e. solving the bond-based PD equations \eref{PDeqn} from time $T_t$ to time $T$ and we report the predicted macro-scale displacement of the middle unit cell in Figure~\ref{fig:test_mid_disp_1D}.

Due to the accumulation of error, the prediction error of the micro-modulus without the energy constraint becomes more pronounced whereas the results corresponding to the use of the energy constraint still maintain high accuracy for all time. 

\begin{figure}[htbp]
\centering
\scalebox{0.6}{\begin{tikzpicture}[]
\begin{axis}[
  height = {101.6mm},
  ylabel = {test error},
  xmin = {0.00015580000000000002},
  xmax = {0.00030419999999999997},
  ymax = {0.6589969806870434},
  ymode = {log},
  xlabel = {$T_t \textrm{ (s})$},
  unbounded coords=jump,scaled x ticks = false,xlabel style = {font = {\fontsize{11 pt}{14.3 pt}\selectfont}, color = {rgb,1:red,0.00000000;green,0.00000000;blue,0.00000000}, draw opacity = 1.0, rotate = 0.0},xmajorgrids = true,xtick = {0.000175,0.0002,0.00022500000000000002,0.00025,0.000275,0.00030000000000000003},xticklabels = {$0.000175$,$0.000200$,$0.000225$,$0.000250$,$0.000275$,$0.000300$},xtick align = inside,xticklabel style = {font = {\fontsize{8 pt}{10.4 pt}\selectfont}, color = {rgb,1:red,0.00000000;green,0.00000000;blue,0.00000000}, draw opacity = 1.0, rotate = 0.0},x grid style = {color = {rgb,1:red,0.00000000;green,0.00000000;blue,0.00000000},
draw opacity = 0.1,
line width = 0.5,
solid},axis x line* = left,x axis line style = {color = {rgb,1:red,0.00000000;green,0.00000000;blue,0.00000000},
draw opacity = 1.0,
line width = 1,
solid},scaled y ticks = false,ylabel style = {font = {\fontsize{11 pt}{14.3 pt}\selectfont}, color = {rgb,1:red,0.00000000;green,0.00000000;blue,0.00000000}, draw opacity = 1.0, rotate = 0.0},log basis y=10,ymajorgrids = true,ytick = {0.001,0.0031622776601683794,0.01,0.03162277660168379,0.1,0.31622776601683794},yticklabels = {$10^{-3.0}$,$10^{-2.5}$,$10^{-2.0}$,$10^{-1.5}$,$10^{-1.0}$,$10^{-0.5}$},ytick align = inside,yticklabel style = {font = {\fontsize{8 pt}{10.4 pt}\selectfont}, color = {rgb,1:red,0.00000000;green,0.00000000;blue,0.00000000}, draw opacity = 1.0, rotate = 0.0},y grid style = {color = {rgb,1:red,0.00000000;green,0.00000000;blue,0.00000000},
draw opacity = 0.1,
line width = 0.5,
solid},axis y line* = left,y axis line style = {color = {rgb,1:red,0.00000000;green,0.00000000;blue,0.00000000},
draw opacity = 1.0,
line width = 1,
solid},    xshift = 0.0mm,
    yshift = 0.0mm,
    axis background/.style={fill={rgb,1:red,1.00000000;green,1.00000000;blue,1.00000000}}
,legend style = {color = {rgb,1:red,0.00000000;green,0.00000000;blue,0.00000000},
draw opacity = 1.0,
line width = 1,
solid,fill = {rgb,1:red,1.00000000;green,1.00000000;blue,1.00000000},fill opacity = 1.0,text opacity = 1.0,font = {\fontsize{8 pt}{10.4 pt}\selectfont}},colorbar style={title=},
  ymin = {0.0004689468182146246},
  width = {152.4mm}
]

\addplot+[
  color = {rgb,1:red,0.00000000;green,0.60560316;blue,0.97868012},
draw opacity = 1.0,
line width = 1,
solid,mark = none,
mark size = 2.0,
mark options = {
            color = {rgb,1:red,0.00000000;green,0.00000000;blue,0.00000000}, draw opacity = 1.0,
            fill = {rgb,1:red,0.00000000;green,0.60560316;blue,0.97868012}, fill opacity = 1.0,
            line width = 1,
            rotate = 0,
            solid
        }
] coordinates {
  (0.00016, 0.5367794517362223)
  (0.00017, 0.06557961181680992)
  (0.00017999999999999998, 0.045597436149506024)
  (0.00019, 0.007901166188006403)
  (0.0002, 0.004750491227137314)
  (0.00021, 0.0023596453845638314)
  (0.0003, 0.0006192887060355528)
};

\addplot+[
  color = {rgb,1:red,0.88887350;green,0.43564919;blue,0.27812294},
draw opacity = 1.0,
line width = 1,
solid,mark = none,
mark size = 2.0,
mark options = {
            color = {rgb,1:red,0.00000000;green,0.00000000;blue,0.00000000}, draw opacity = 1.0,
            fill = {rgb,1:red,0.88887350;green,0.43564919;blue,0.27812294}, fill opacity = 1.0,
            line width = 1,
            rotate = 0,
            solid
        }
] coordinates {
  (0.00016, 0.014631594504956085)
  (0.00017, 0.008212360595867085)
  (0.00017999999999999998, 0.0012196296567273926)
  (0.00019, 0.0014768150282089118)
  (0.0002, 0.0007083124894100021)
  (0.00021, 0.0005757197603348187)
  (0.0003, 0.0005791119377546962)
};

\legend{{}{without constraint}, {}{with constraint}}
\end{axis}

\end{tikzpicture}}
\caption{Relative testing error of $\ddot{\bar u}$ for different size of training data}
\label{fig:test_error_1D}
\end{figure}
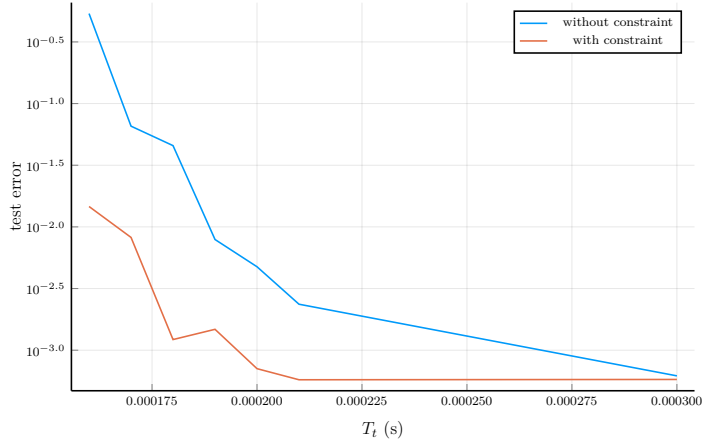

In order to better show the improvements of the energy constraint on the performance of our algorithm, we vary the size of the training data (i.e. vary the value of $T_t$) and compute the relative prediction error in $\ell^2$ norm.
The results shown in Figure~\ref{fig:test_error_1D} indicate that, when the amount of training data is small, the energy constraint provides more accurate predictions, whereas when the amount of training data increases, the performance of the learning algorithm without the energy constraint is similar to that with the energy constraint.

\begin{figure}[htbp]
\centering
\scalebox{0.65}{\input{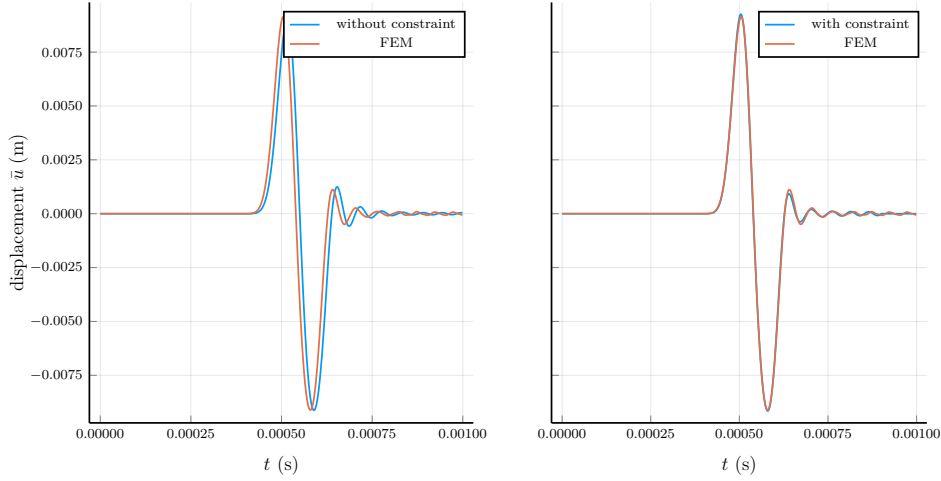}}
\caption{Validation results on the displacement of the middle unit cell}
\label{fig:valid_mid_disp_1D}
\end{figure}

We further validate the performance of the learning algorithm by using the optimal micro-modulus function learned in the previous experiments to predict the deformation of the rod under different loads. We change the time-dependent displacement boundary condition to $u_{bc}(t) = u_0 a_0 \sin(2\pi t/T_s)[1-H(t-T_s)]$, solve the bond-based PD equations with the optimal micro-modulus, and compare the displacement solutions with the macro-scale displacement calculated using FEM. The comparison results are shown in Figure~\ref{fig:valid_mid_disp_1D}: we observe that the bond-based PD model successfully describes the deformation of the periodic heterogeneous rod under a different loads. This shows that our algorithm generalizes well for different loading scenarios.

\subsection{Two-dimensional test cases}

We consider the two-dimensional plane stress problem described in {\S\ref{chap:2Delastic}} with geometry parameters $L = 1 \textrm{m}$ and $l_2 = 1/3 \textrm{m}$ and material properties values are: $E_s = 200 \textrm{GPa}$, $E_c = 5 \textrm{GPa}$, and $\rho = 8000 \textrm{kg/m}^3$. The time-dependent displacement boundary condition used for FEMsimulations is given by
\begin{align}
    \bar u_x(L,y) = u_0 a_0 t^6(t-T_s)^6[1-H(t-T_s)]  \qquad \quad 0\leq y\leq L,
    \label{eqn:2DBC}
\end{align}
where $u_0 = 1\times10^{-3} \mathrm{m}$, $a_0$ is a scaling factor, $H$ is the Heaviside function and $T_s = 0.0785 \textrm{ms}$. 
The macro-scale data set is generated the same way as in the previous section with $T=5\times10^{-4}\rm s$ and $T_t=0.8\times10^{-4}\rm s$.
We choose the horizon of the PD model to be $\epsilon = 6l$ and, as done in the one-dimensional case, we do not consider macro-scale data of points whose nonlocal neighborhood is not contained in the domain. 

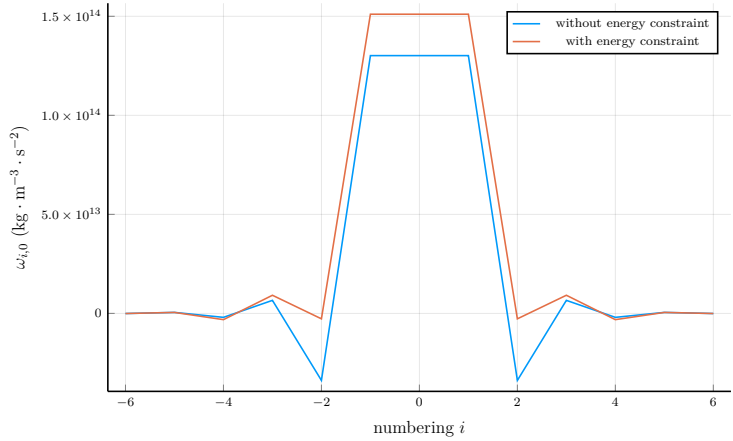
\begin{figure}[htbp]
\centering
\scalebox{0.6}{\begin{tikzpicture}[]
\begin{axis}[
  height = {101.6mm},
  ylabel = {$\omega_{i,0}\textrm{ }(\textrm{kg}\cdot \textrm{m}^{-3}\cdot\textrm{s}^{-2})$\\\color{white}.\\\color{white}.\\\color{white}.},
  xmin = {-6.36},
  xmax = {6.36},
  ymax = {1.565961167615651e14},
  xlabel = {$\textrm{numbering }i$},
  unbounded coords=jump,scaled x ticks = false,xlabel style = {font = {\fontsize{11 pt}{14.3 pt}\selectfont}, color = {rgb,1:red,0.00000000;green,0.00000000;blue,0.00000000}, draw opacity = 1.0, rotate = 0.0},xmajorgrids = true,xtick = {-6.0,-4.0,-2.0,0.0,2.0,4.0,6.0},xticklabels = {$-6$,$-4$,$-2$,$0$,$2$,$4$,$6$},xtick align = inside,xticklabel style = {font = {\fontsize{8 pt}{10.4 pt}\selectfont}, color = {rgb,1:red,0.00000000;green,0.00000000;blue,0.00000000}, draw opacity = 1.0, rotate = 0.0},x grid style = {color = {rgb,1:red,0.00000000;green,0.00000000;blue,0.00000000},
draw opacity = 0.1,
line width = 0.5,
solid},axis x line* = left,x axis line style = {color = {rgb,1:red,0.00000000;green,0.00000000;blue,0.00000000},
draw opacity = 1.0,
line width = 1,
solid},scaled y ticks = false,ylabel style = {align=center,font = {\fontsize{11 pt}{14.3 pt}\selectfont}, color = {rgb,1:red,0.00000000;green,0.00000000;blue,0.00000000}, draw opacity = 1.0, rotate = 0.0},ymajorgrids = true,ytick = {0.0,5.0e13,1.0e14,1.5e14},yticklabels = {$0$,$5.0\times10^{13}$,$1.0\times10^{14}$,$1.5\times10^{14}$},ytick align = inside,yticklabel style = {font = {\fontsize{8 pt}{10.4 pt}\selectfont}, color = {rgb,1:red,0.00000000;green,0.00000000;blue,0.00000000}, draw opacity = 1.0, rotate = 0.0},y grid style = {color = {rgb,1:red,0.00000000;green,0.00000000;blue,0.00000000},
draw opacity = 0.1,
line width = 0.5,
solid},axis y line* = left,y axis line style = {color = {rgb,1:red,0.00000000;green,0.00000000;blue,0.00000000},
draw opacity = 1.0,
line width = 1,
solid},    xshift = 0.0mm,
    yshift = 0.0mm,
    axis background/.style={fill={rgb,1:red,1.00000000;green,1.00000000;blue,1.00000000}}
,legend style = {color = {rgb,1:red,0.00000000;green,0.00000000;blue,0.00000000},
draw opacity = 1.0,
line width = 1,
solid,fill = {rgb,1:red,1.00000000;green,1.00000000;blue,1.00000000},fill opacity = 1.0,text opacity = 1.0,font = {\fontsize{8 pt}{10.4 pt}\selectfont}},colorbar style={title=},
  ymin = {-3.93767588615514e13},
  width = {152.40000000000003mm}
]

\addplot+[
  color = {rgb,1:red,0.00000000;green,0.60560316;blue,0.97868012},
draw opacity = 1.0,
line width = 1,
solid,mark = none,
mark size = 2.0,
mark options = {
            color = {rgb,1:red,0.00000000;green,0.00000000;blue,0.00000000}, draw opacity = 1.0,
            fill = {rgb,1:red,0.00000000;green,0.60560316;blue,0.97868012}, fill opacity = 1.0,
            line width = 1,
            rotate = 0,
            solid
        }
] coordinates {
  (-6.0, -7.513822515192159e10)
  (-5.0, 4.9739533832243774e11)
  (-4.0, -2.0039241762260928e12)
  (-3.0, 6.574291421884126e12)
  (-2.0, -3.383035672127452e13)
  (-1.0, 1.3014365161786655e14)
  (0.0, 1.3014365161786655e14)
  (1.0, 1.3014365161786655e14)
  (2.0, -3.383035672127452e13)
  (3.0, 6.574291421884126e12)
  (4.0, -2.0039241762260928e12)
  (5.0, 4.9739533832243774e11)
  (6.0, -7.513822515192159e10)
};

\addplot+[
  color = {rgb,1:red,0.88887350;green,0.43564919;blue,0.27812294},
draw opacity = 1.0,
line width = 1,
solid,mark = none,
mark size = 2.0,
mark options = {
            color = {rgb,1:red,0.00000000;green,0.00000000;blue,0.00000000}, draw opacity = 1.0,
            fill = {rgb,1:red,0.88887350;green,0.43564919;blue,0.27812294}, fill opacity = 1.0,
            line width = 1,
            rotate = 0,
            solid
        }
] coordinates {
  (-6.0, -7.052040722628088e10)
  (-5.0, 4.8382705704811945e11)
  (-4.0, -3.1628803572515615e12)
  (-3.0, 9.151536422520934e12)
  (-2.0, -2.7482111000628433e12)
  (-1.0, 1.5104971462128822e14)
  (0.0, 1.5104971462128822e14)
  (1.0, 1.5104971462128822e14)
  (2.0, -2.7482111000628433e12)
  (3.0, 9.151536422520934e12)
  (4.0, -3.1628803572515615e12)
  (5.0, 4.8382705704811945e11)
  (6.0, -7.052040722628088e10)
};

\legend{{}{without energy constraint}, {}{with energy constraint}}
\end{axis}

\end{tikzpicture}}
\caption{Discrete micro-modulus $\omega_{i,0}$ for two-dimensional elasticity}
\label{fig:micro-modulus_2D}
\end{figure}

\begin{figure}[htbp]
\centering
\scalebox{0.65}{\input{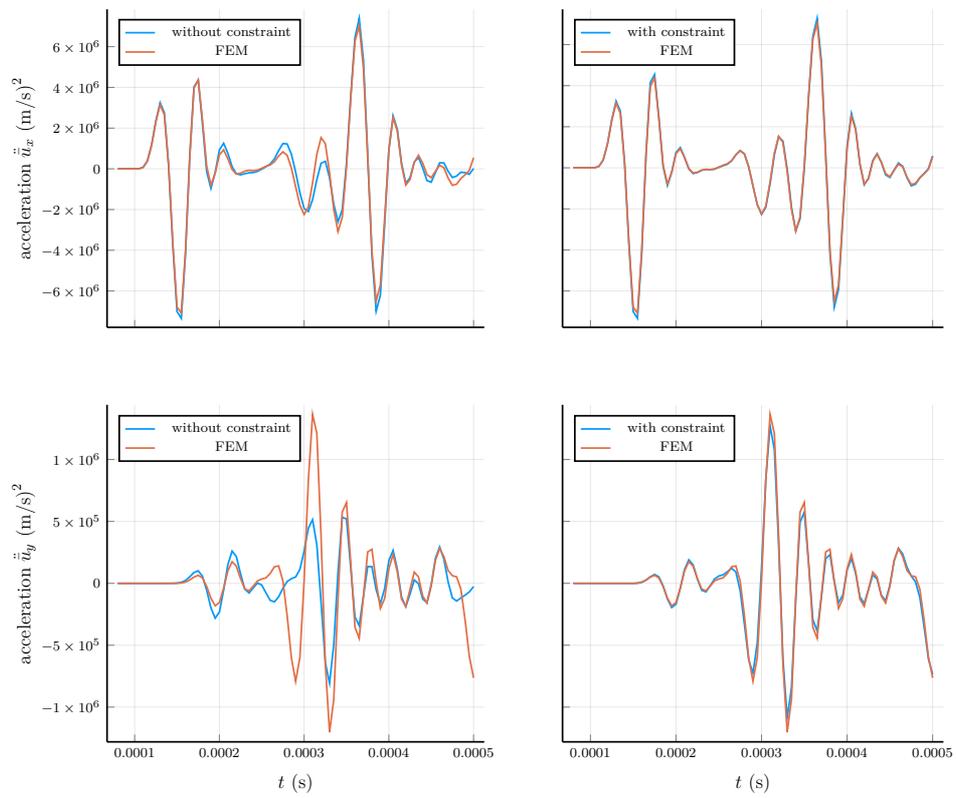}}
\caption{Testing results on the acceleration of the middle unit cell}
\label{fig:test_mid_acc_2D}
\end{figure}

\begin{figure}[htbp]
\centering
\scalebox{0.65}{\input{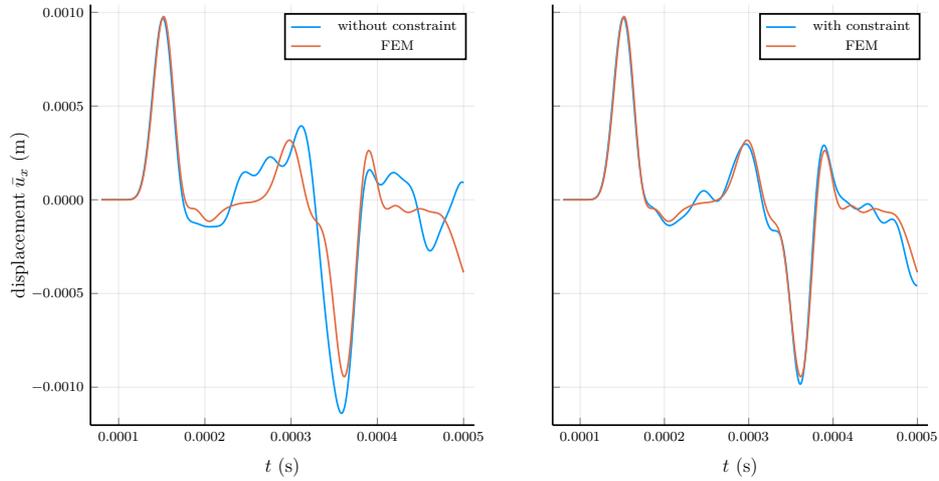}}
\caption{Testing results on the displacement of the middle unit cell}
\label{fig:test_mid_disp_2D}
\end{figure}

We select an appropriate weighted $\ell^2$ norm for the objective functions, perform linear regression with and without energy constraint and report the optimal values of the discrete micro-modulus function $\{\omega^*_{i,j}\}_{i,j=-6}^{i,j=6}$ at $j=0$ in Figure~\ref{fig:micro-modulus_2D}. 
Also in this case these two micro-moduli generate positive definite matrices when discretizing equation \eref{PDeqn}.
We test the two optimal micro-moduli on the testing data set and report the predicted macro-scale acceleration $\ddot{\bar u}_{x}$ and $\ddot{\bar u}_{y}$ of the middle unit cell as a function of time in Figure~\ref{fig:test_mid_acc_2D}, compared with corresponding macro-scale FEM solutions.
It is evident that the accuracy of the energy-constrained prediction is better when compared to the unconstrained one, especially for the $y$-component of acceleration. This is not unexpected since the magnitude of the training data is higher in the $x$ direction. As reported in the previous section, we also solve the bond-based PD equations using the two optimal micro-moduli in $(T_t,T)$ and report the corresponding ${\bar u}_x$ together with the reference FEM solutions in Figure~\ref{fig:test_mid_disp_2D}. In Figures~\ref{fig:test_errorx_2D} and \ref{fig:test_errory_2D} we report the relative testing error (defined as in the previous section) in correspondence of different sizes of the training data set (i.e. different values of $T_t$) for $\ddot{\bar u}_x$ and $\ddot{\bar u}_y$ respectively. Similar considerations as in the one-dimensional case can be inferred.
\begin{figure}[htbp]
\centering
\scalebox{0.6}{\begin{tikzpicture}[]
\begin{axis}[
  height = {101.6mm},
  ylabel = {test error},
  xmin = {6.76e-5},
  xmax = {0.0001524},
  ymax = {1.0929785883937877},
  ymode = {log},
  xlabel = {$T_t \textrm{ (s)}$},
  unbounded coords=jump,scaled x ticks = false,xlabel style = {font = {\fontsize{11 pt}{14.3 pt}\selectfont}, color = {rgb,1:red,0.00000000;green,0.00000000;blue,0.00000000}, draw opacity = 1.0, rotate = 0.0},xmajorgrids = true,xtick = {8.0e-5,0.0001,0.00012000000000000002,0.00014000000000000001},xticklabels = {$8.0\times10^{-5}$,$1.0\times10^{-4}$,$1.2\times10^{-4}$,$1.4\times10^{-4}$},xtick align = inside,xticklabel style = {font = {\fontsize{8 pt}{10.4 pt}\selectfont}, color = {rgb,1:red,0.00000000;green,0.00000000;blue,0.00000000}, draw opacity = 1.0, rotate = 0.0},x grid style = {color = {rgb,1:red,0.00000000;green,0.00000000;blue,0.00000000},
draw opacity = 0.1,
line width = 0.5,
solid},axis x line* = left,x axis line style = {color = {rgb,1:red,0.00000000;green,0.00000000;blue,0.00000000},
draw opacity = 1.0,
line width = 1,
solid},scaled y ticks = false,ylabel style = {font = {\fontsize{11 pt}{14.3 pt}\selectfont}, color = {rgb,1:red,0.00000000;green,0.00000000;blue,0.00000000}, draw opacity = 1.0, rotate = 0.0},log basis y=10,ymajorgrids = true,ytick = {0.01,0.03162277660168379,0.1,0.31622776601683794,1.0},yticklabels = {$10^{-2.0}$,$10^{-1.5}$,$10^{-1.0}$,$10^{-0.5}$,$10^{0.0}$},ytick align = inside,yticklabel style = {font = {\fontsize{8 pt}{10.4 pt}\selectfont}, color = {rgb,1:red,0.00000000;green,0.00000000;blue,0.00000000}, draw opacity = 1.0, rotate = 0.0},y grid style = {color = {rgb,1:red,0.00000000;green,0.00000000;blue,0.00000000},
draw opacity = 0.1,
line width = 0.5,
solid},axis y line* = left,y axis line style = {color = {rgb,1:red,0.00000000;green,0.00000000;blue,0.00000000},
draw opacity = 1.0,
line width = 1,
solid},    xshift = 0.0mm,
    yshift = 0.0mm,
    axis background/.style={fill={rgb,1:red,1.00000000;green,1.00000000;blue,1.00000000}}
,legend style = {color = {rgb,1:red,0.00000000;green,0.00000000;blue,0.00000000},
draw opacity = 1.0,
line width = 1,
solid,fill = {rgb,1:red,1.00000000;green,1.00000000;blue,1.00000000},fill opacity = 1.0,text opacity = 1.0,font = {\fontsize{8 pt}{10.4 pt}\selectfont}},colorbar style={title=},
  ymin = {0.007168103714513083},
  width = {152.4mm}
]

\addplot+[
  color = {rgb,1:red,0.00000000;green,0.60560316;blue,0.97868012},
draw opacity = 1.0,
line width = 1,
solid,mark = none,
mark size = 2.0,
mark options = {
            color = {rgb,1:red,0.00000000;green,0.00000000;blue,0.00000000}, draw opacity = 1.0,
            fill = {rgb,1:red,0.00000000;green,0.60560316;blue,0.97868012}, fill opacity = 1.0,
            line width = 1,
            rotate = 0,
            solid
        }
] coordinates {
  (7.000000000000001e-5, 0.9480314965912725)
  (7.999999999999999e-5, 0.11459323977195686)
  (8.999999999999999e-5, 0.029125022202065372)
  (9.999999999999999e-5, 0.015500436413710858)
  (0.00010999999999999999, 0.021188512970049807)
  (0.00015, 0.008264054419624967)
};

\addplot+[
  color = {rgb,1:red,0.88887350;green,0.43564919;blue,0.27812294},
draw opacity = 1.0,
line width = 1,
solid,mark = none,
mark size = 2.0,
mark options = {
            color = {rgb,1:red,0.00000000;green,0.00000000;blue,0.00000000}, draw opacity = 1.0,
            fill = {rgb,1:red,0.88887350;green,0.43564919;blue,0.27812294}, fill opacity = 1.0,
            line width = 1,
            rotate = 0,
            solid
        }
] coordinates {
  (7.000000000000001e-5, 0.1472562617329486)
  (7.999999999999999e-5, 0.036488044028037)
  (8.999999999999999e-5, 0.01720047200270083)
  (9.999999999999999e-5, 0.015342140827378423)
  (0.00010999999999999999, 0.018873825977195422)
  (0.00015, 0.008711460572655132)
};

\legend{{}{without constraint}, {}{with constraint}}
\end{axis}

\end{tikzpicture}}
\caption{Relative testing error of $\ddot{\bar u}_x$ for different size of training data}
\label{fig:test_errorx_2D}
\end{figure}
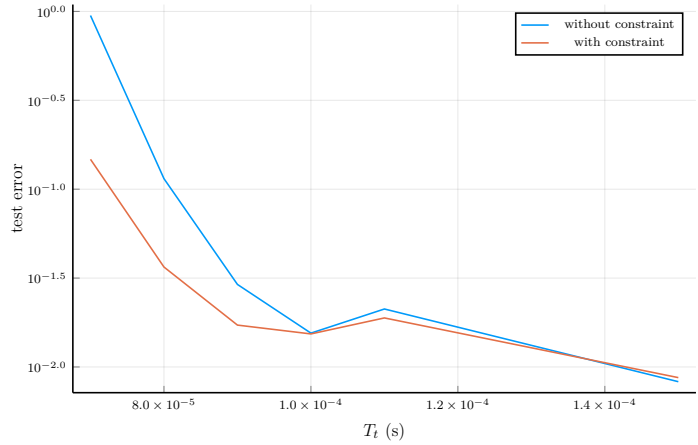
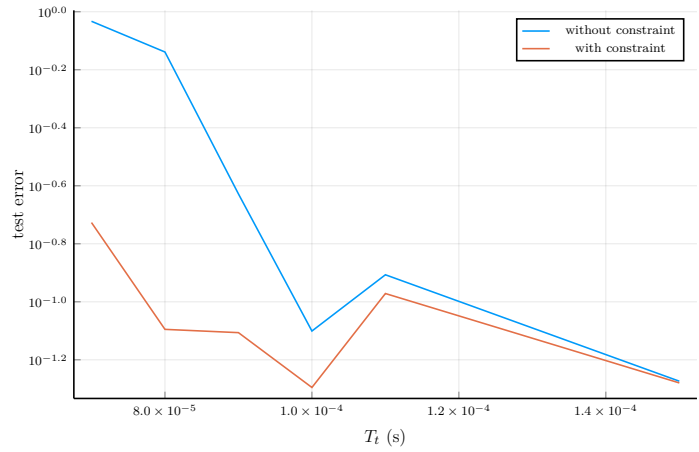
\begin{figure}[htbp]
\centering
\scalebox{0.6}{\begin{tikzpicture}[]
\begin{axis}[
  height = {101.6mm},
  ylabel = {test error},
  xmin = {6.76e-5},
  xmax = {0.0001524},
  ymax = {1.011796933587817},
  ymode = {log},
  xlabel = {$T_t \textrm{ (s)}$},
  unbounded coords=jump,scaled x ticks = false,xlabel style = {font = {\fontsize{11 pt}{14.3 pt}\selectfont}, color = {rgb,1:red,0.00000000;green,0.00000000;blue,0.00000000}, draw opacity = 1.0, rotate = 0.0},xmajorgrids = true,xtick = {8.0e-5,0.0001,0.00012000000000000002,0.00014000000000000001},xticklabels = {$8.0\times10^{-5}$,$1.0\times10^{-4}$,$1.2\times10^{-4}$,$1.4\times10^{-4}$},xtick align = inside,xticklabel style = {font = {\fontsize{8 pt}{10.4 pt}\selectfont}, color = {rgb,1:red,0.00000000;green,0.00000000;blue,0.00000000}, draw opacity = 1.0, rotate = 0.0},x grid style = {color = {rgb,1:red,0.00000000;green,0.00000000;blue,0.00000000},
draw opacity = 0.1,
line width = 0.5,
solid},axis x line* = left,x axis line style = {color = {rgb,1:red,0.00000000;green,0.00000000;blue,0.00000000},
draw opacity = 1.0,
line width = 1,
solid},scaled y ticks = false,ylabel style = {font = {\fontsize{11 pt}{14.3 pt}\selectfont}, color = {rgb,1:red,0.00000000;green,0.00000000;blue,0.00000000}, draw opacity = 1.0, rotate = 0.0},log basis y=10,ymajorgrids = true,ytick = {0.0630957344480193,0.1,0.15848931924611134,0.25118864315095796,0.3981071705534972,0.6309573444801932,1.0},yticklabels = {$10^{-1.2}$,$10^{-1.0}$,$10^{-0.8}$,$10^{-0.6}$,$10^{-0.4}$,$10^{-0.2}$,$10^{0.0}$},ytick align = inside,yticklabel style = {font = {\fontsize{8 pt}{10.4 pt}\selectfont}, color = {rgb,1:red,0.00000000;green,0.00000000;blue,0.00000000}, draw opacity = 1.0, rotate = 0.0},y grid style = {color = {rgb,1:red,0.00000000;green,0.00000000;blue,0.00000000},
draw opacity = 0.1,
line width = 0.5,
solid},axis y line* = left,y axis line style = {color = {rgb,1:red,0.00000000;green,0.00000000;blue,0.00000000},
draw opacity = 1.0,
line width = 1,
solid},    xshift = 0.0mm,
    yshift = 0.0mm,
    axis background/.style={fill={rgb,1:red,1.00000000;green,1.00000000;blue,1.00000000}}
,legend style = {color = {rgb,1:red,0.00000000;green,0.00000000;blue,0.00000000},
draw opacity = 1.0,
line width = 1,
solid,fill = {rgb,1:red,1.00000000;green,1.00000000;blue,1.00000000},fill opacity = 1.0,text opacity = 1.0,font = {\fontsize{8 pt}{10.4 pt}\selectfont}},colorbar style={title=},
  ymin = {0.04643031889945853},
  width = {152.4mm}
]

\addplot+[
  color = {rgb,1:red,0.00000000;green,0.60560316;blue,0.97868012},
draw opacity = 1.0,
line width = 1,
solid,mark = none,
mark size = 2.0,
mark options = {
            color = {rgb,1:red,0.00000000;green,0.00000000;blue,0.00000000}, draw opacity = 1.0,
            fill = {rgb,1:red,0.00000000;green,0.60560316;blue,0.97868012}, fill opacity = 1.0,
            line width = 1,
            rotate = 0,
            solid
        }
] coordinates {
  (7.000000000000001e-5, 0.9272934217713802)
  (7.999999999999999e-5, 0.726671371452769)
  (8.999999999999999e-5, 0.23551377425940642)
  (9.999999999999999e-5, 0.07931287003679556)
  (0.00010999999999999999, 0.12395467398409199)
  (0.00015, 0.05324573378877649)
};

\addplot+[
  color = {rgb,1:red,0.88887350;green,0.43564919;blue,0.27812294},
draw opacity = 1.0,
line width = 1,
solid,mark = none,
mark size = 2.0,
mark options = {
            color = {rgb,1:red,0.00000000;green,0.00000000;blue,0.00000000}, draw opacity = 1.0,
            fill = {rgb,1:red,0.88887350;green,0.43564919;blue,0.27812294}, fill opacity = 1.0,
            line width = 1,
            rotate = 0,
            solid
        }
] coordinates {
  (7.000000000000001e-5, 0.18763081500047074)
  (7.999999999999999e-5, 0.08035964918031636)
  (8.999999999999999e-5, 0.07829742926121745)
  (9.999999999999999e-5, 0.0506614769230605)
  (0.00010999999999999999, 0.10679723923831419)
  (0.00015, 0.05252443172597384)
};

\legend{{}{without constraint}, {}{with constraint}}
\end{axis}

\end{tikzpicture}}
\caption{Relative testing error of $\ddot{\bar u}_y$ for different size of training data}
\label{fig:test_errory_2D}
\end{figure}
\begin{figure}[htbp]
\centering
\scalebox{0.65}{\input{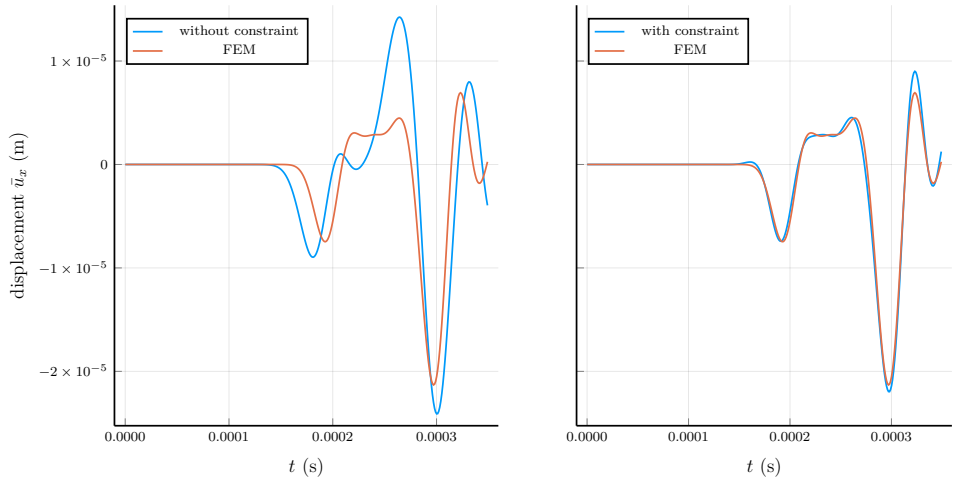}}
\caption{Validation results on $x$ displacement of the middle unit cell}
\label{fig:valid_mid_dispx_2D}
\end{figure}

Finally, we test the optimal micro-modulus functions used in the previous tests on another macro-scale data set generated using a different load. Specifically, we generate the training and validation data set with the following boundary condition, which corresponds to a shear load
\begin{align}
    \bar u_y(L,y) = u_0 a_0 t^6(t-T_s)^6[1-H(t-T_s)]  \qquad \quad 0\leq y\leq L. \notag 
\end{align}
The prediction results for the $x$-displacement are shown in Figure~\ref{fig:valid_mid_dispx_2D}; once again, the accuracy of the predictions implies that our energy-constrained learning algorithm generalizes well to loads that are different from the one used for training.
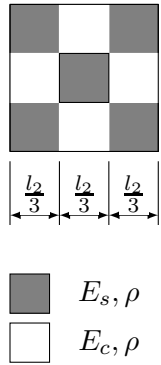
\begin{figure}[htbp]
\centering
\begin{tikzpicture}[scale=0.65]
    \tikzstyle{pt}=[circle, fill=black, inner sep=0pt, minimum size=5pt]
    
    
    
    
    \filldraw[fill=black!50, draw=black] (4,3) rectangle (7,0);
    \filldraw[fill=white!50, draw=white] (4,1) rectangle (7,2);
    \filldraw[fill=white!50, draw=white] (5,3) rectangle (6,0);
    \filldraw[fill=black!50, draw=black] (5,2) rectangle (6,1);
    \draw (4,0) -- (4,3);
    \draw (4,0) -- (7,0);
    \draw (7,0) -- (7,3);
    \draw (4,3) -- (7,3);
    
    
    
    
    
    
    
    \draw[latex-latex] (4,-1.3) -- (5,-1.3) node[midway , above,  fill=white, inner sep=1pt] {$\frac{l_2}{3}$};
    \draw (4,-1.5) -- (4,-0.2);
    \draw (5,-1.5) -- (5,-0.2);
    
    \draw[latex-latex] (6,-1.3) -- (7,-1.3) node[midway , above,  fill=white, inner sep=1pt] {$\frac{l_2}{3}$};
    \draw (6,-1.5) -- (6,-0.2);
    \draw (7,-1.5) -- (7,-0.2);
    
    \draw[latex-latex] (5,-1.3) -- (6,-1.3) node[midway , above,  fill=white, inner sep=1pt] {$\frac{l_2}{3}$};
    
    \filldraw[fill=black!50, draw=black] (4,-1-1.5) rectangle (4.8,-1.8-1.5);
    \draw (6,-1.4-1.5) node() {$ E_s,\rho$};
    
    \filldraw[fill=white, draw=black] (4,-2-1.5) rectangle (4.8,-2.8-1.5);
    \draw (6,-2.4-1.5) node() {$ E_c,\rho$};
    

    

\end{tikzpicture}
\caption{Two-dimensional heterogeneous unit cell}
\label{fig:newmicrostructure}
\end{figure}
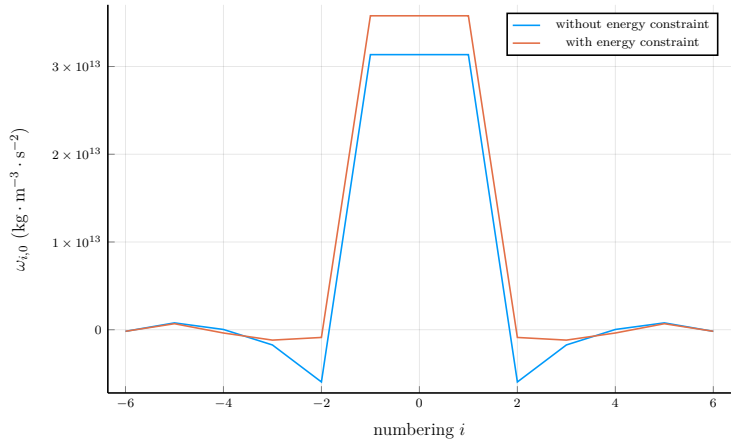
\begin{figure}[htbp]
\centering
\scalebox{0.6}{\begin{tikzpicture}[]
\begin{axis}[
  height = {101.6mm},
  ylabel = {$\omega_{i,0}\textrm{ }(\textrm{kg}\cdot \textrm{m}^{-3}\cdot\textrm{s}^{-2})$\\\color{white}.\\\color{white}.\\\color{white}.},
  xmin = {-6.36},
  xmax = {6.36},
  ymax = {3.7032933364669305e13},
  xlabel = {$\textrm{numbering }i$},
  unbounded coords=jump,scaled x ticks = false,xlabel style = {font = {\fontsize{11 pt}{14.3 pt}\selectfont}, color = {rgb,1:red,0.00000000;green,0.00000000;blue,0.00000000}, draw opacity = 1.0, rotate = 0.0},xmajorgrids = true,xtick = {-6.0,-4.0,-2.0,0.0,2.0,4.0,6.0},xticklabels = {$-6$,$-4$,$-2$,$0$,$2$,$4$,$6$},xtick align = inside,xticklabel style = {font = {\fontsize{8 pt}{10.4 pt}\selectfont}, color = {rgb,1:red,0.00000000;green,0.00000000;blue,0.00000000}, draw opacity = 1.0, rotate = 0.0},x grid style = {color = {rgb,1:red,0.00000000;green,0.00000000;blue,0.00000000},
draw opacity = 0.1,
line width = 0.5,
solid},axis x line* = left,x axis line style = {color = {rgb,1:red,0.00000000;green,0.00000000;blue,0.00000000},
draw opacity = 1.0,
line width = 1,
solid},scaled y ticks = false,ylabel style = {align=center,font = {\fontsize{11 pt}{14.3 pt}\selectfont}, color = {rgb,1:red,0.00000000;green,0.00000000;blue,0.00000000}, draw opacity = 1.0, rotate = 0.0},ymajorgrids = true,ytick = {0.0,1.0e13,2.0e13,3.0e13},yticklabels = {$0$,$1\times10^{13}$,$2\times10^{13}$,$3\times10^{13}$},ytick align = inside,yticklabel style = {font = {\fontsize{8 pt}{10.4 pt}\selectfont}, color = {rgb,1:red,0.00000000;green,0.00000000;blue,0.00000000}, draw opacity = 1.0, rotate = 0.0},y grid style = {color = {rgb,1:red,0.00000000;green,0.00000000;blue,0.00000000},
draw opacity = 0.1,
line width = 0.5,
solid},axis y line* = left,y axis line style = {color = {rgb,1:red,0.00000000;green,0.00000000;blue,0.00000000},
draw opacity = 1.0,
line width = 1,
solid},    xshift = 0.0mm,
    yshift = 0.0mm,
    axis background/.style={fill={rgb,1:red,1.00000000;green,1.00000000;blue,1.00000000}}
,legend style = {color = {rgb,1:red,0.00000000;green,0.00000000;blue,0.00000000},
draw opacity = 1.0,
line width = 1,
solid,fill = {rgb,1:red,1.00000000;green,1.00000000;blue,1.00000000},fill opacity = 1.0,text opacity = 1.0,font = {\fontsize{8 pt}{10.4 pt}\selectfont}},colorbar style={title=},
  ymin = {-7.205212329662865e12},
  width = {152.4mm}
]

\addplot+[
  color = {rgb,1:red,0.00000000;green,0.60560316;blue,0.97868012},
draw opacity = 1.0,
line width = 1,
solid,mark = none,
mark size = 2.0,
mark options = {
            color = {rgb,1:red,0.00000000;green,0.00000000;blue,0.00000000}, draw opacity = 1.0,
            fill = {rgb,1:red,0.00000000;green,0.60560316;blue,0.97868012}, fill opacity = 1.0,
            line width = 1,
            rotate = 0,
            solid
        }
] coordinates {
  (-6.0, -1.8312110470257907e11)
  (-5.0, 7.815558172224883e11)
  (-4.0, 3.0259704810358135e10)
  (-3.0, -1.735159862436489e12)
  (-2.0, -5.953189338313842e12)
  (-1.0, 3.1353329868921793e13)
  (0.0, 3.1353329868921793e13)
  (1.0, 3.1353329868921793e13)
  (2.0, -5.953189338313842e12)
  (3.0, -1.735159862436489e12)
  (4.0, 3.0259704810358135e10)
  (5.0, 7.815558172224883e11)
  (6.0, -1.8312110470257907e11)
};

\addplot+[
  color = {rgb,1:red,0.88887350;green,0.43564919;blue,0.27812294},
draw opacity = 1.0,
line width = 1,
solid,mark = none,
mark size = 2.0,
mark options = {
            color = {rgb,1:red,0.00000000;green,0.00000000;blue,0.00000000}, draw opacity = 1.0,
            fill = {rgb,1:red,0.88887350;green,0.43564919;blue,0.27812294}, fill opacity = 1.0,
            line width = 1,
            rotate = 0,
            solid
        }
] coordinates {
  (-6.0, -1.7843178436780353e11)
  (-5.0, 6.846571635205648e11)
  (-4.0, -3.724649545641827e11)
  (-3.0, -1.1824523812090579e12)
  (-2.0, -8.797114343232277e11)
  (-1.0, 3.578091037332028e13)
  (0.0, 3.578091037332028e13)
  (1.0, 3.578091037332028e13)
  (2.0, -8.797114343232277e11)
  (3.0, -1.1824523812090579e12)
  (4.0, -3.724649545641827e11)
  (5.0, 6.846571635205648e11)
  (6.0, -1.7843178436780353e11)
};

\legend{{}{without energy constraint}, {}{with energy constraint}}
\end{axis}

\end{tikzpicture}}
\caption{Discrete micro-modulus $\omega_{i,0}$ for two-dimensional elasticity}
\label{fig:micro-modulus_2D_2}
\end{figure}
\begin{figure}[htbp]
\centering
\scalebox{0.65}{\input{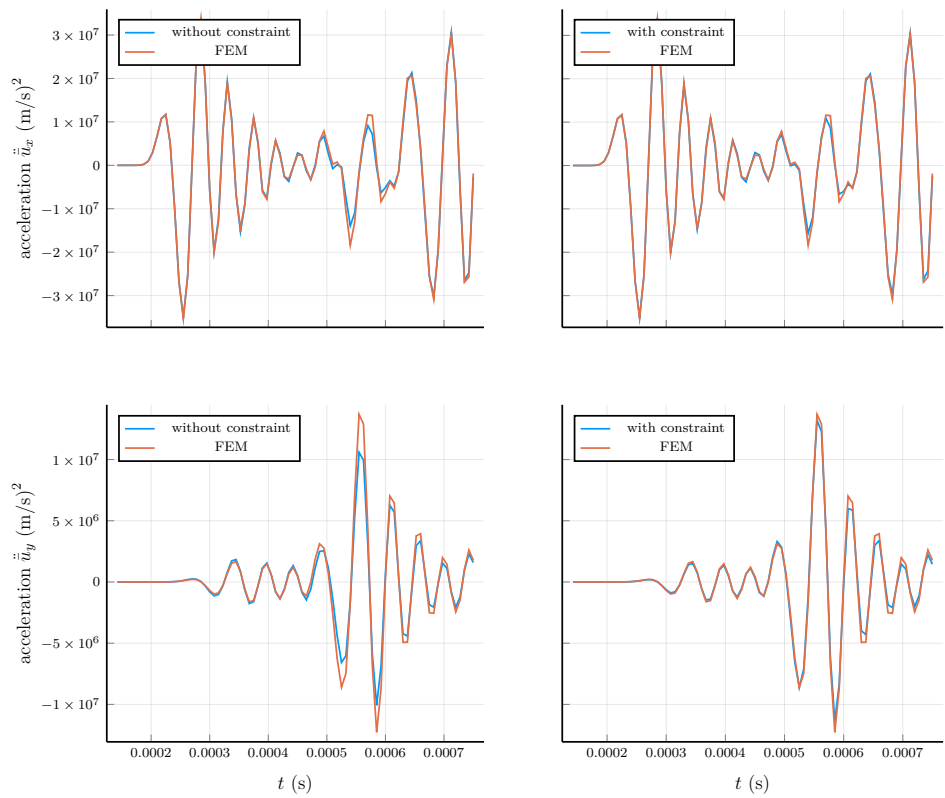}}
\caption{Testing results on the acceleration of the middle unit cell}
\label{fig:test_mid_acc_2D_2}
\end{figure}

With the purpose of proving the generality of our algorithm, we consider different medium of the same size as the one in Figure~\ref{fig:2DIllustration}, but with different micro-structure (see the illustration in Figure~\ref{fig:newmicrostructure}). We use the time-dependent boundary condition reported in \eref{2DBC} and generate the macro-scale data which is divided into training set and testing data set with $T = 7.5 \times 10^{-4}$ and $T_t = 2.9\times 10^{-4}$.
We choose the horizon for this medium to be $\epsilon = 6l_2$, and perform linear regression on the training set with and without the energy constraint. The resulting discrete micro-modulus functions are plotted in Figure~\ref{fig:micro-modulus_2D_2}. We use these two discrete micro-modulus functions to predict the macro-scale acceleration of the testing data set and plot the comparison results of the middle unit cell in Figure~\ref{fig:test_mid_acc_2D_2}. The accuracy of our results indicates that the learning algorithm performs well regardless of the micro-structure.

\section{Conclusions and Future Work}
\label{chap:conclusions}
The importance of the PD influence function has recently gained more attention from the PD research community, but what has been missing in the literature is a systematic way for it to be determined for a given material and application setting. In this work, we use a bond-based PD model to describe the linear elastic deformations for materials with periodic heterogeneity. We used a highly-resolved micro-structural FEM to solve the classical elastodynamic equations for a short time periods and upscaled the solutions to an averaged macro-scale deformation. {We used these solutions as training data in a machine-learning framework to identify} the optimal discrete micro-modulus for the PD model. In the regression algorithm, an energy constraint that represents the average elastic modulus of the micro-structure was added to the objective function. The testing results indicated that the homogenized macro-scale deformation can be predicted by the resulting micro-modulus and the energy constraint helps in constructing the PD model with better accuracy using less data.

In the interest of simplicity, this work only focuses on the application of bond-based PD models and an energy constraint corresponding to isotropic deformation, which reduces the regression algorithm to simple linear regression. Future work should include extensions of the algorithm to more complex state-based PD models and constructing additional mathematically- and physically-justified constraints, e.g.\ constraints that can guarantee the well-posedness of resulting PD model and/or constraints that contain more detailed information of the medium's micro-structure.

\section{Acknowledgment}
This work was partially supported by the Sandia National Laboratories (SNL) Laboratory-directed Research and Development program and by the U.S. Department of Energy, Office of Advanced Scientific Computing Research under the Collaboratory on Mathematics and Physics-Informed Learning Machines for Multiscale and Multiphysics Problems (PhILMs) project. SNL is a multimission laboratory managed and operated by National Technology and Engineering Solutions of Sandia, LLC., a wholly owned subsidiary of Honeywell International, Inc., for the U.S. Department of Energy's National Nuclear Security Administration under contract DE-NA-0003525. This paper (SAND2021-0028) describes objective technical results and analysis. Any subjective views or opinions that might be expressed in this paper do not necessarily represent the views of the U.S. Department of Energy or the United States Government. 

\bibliographystyle{abbrvnat}
\bibliography{mybib}
\end{document}